\documentclass[12pt]{article}
\usepackage[utf8]{inputenc}
\usepackage[T1]{fontenc}
\usepackage{theorem}
\usepackage{dsfont}
\usepackage{mathtools}
\usepackage{multirow}
\usepackage{amssymb}
\usepackage{array}
\usepackage{color}
\usepackage{graphicx}
\usepackage{pgfplots}
\usepackage{tikz}
\usetikzlibrary{positioning}
\usepackage[export]{adjustbox}
\usetikzlibrary{shapes.geometric, arrows}
\usepackage{amsmath}
\usepackage{verbatim} 
\usepackage{rotating}
\usepackage[round]{natbib}
\usepackage[colorlinks=true,linkcolor=black, citecolor=blue, urlcolor=blue]{hyperref}
\usepackage{footnote}
\usepackage{tcolorbox}
\usepackage{bookmark}
\makesavenoteenv{tabular}
\makesavenoteenv{table}

\newtheorem{theorem}{Theorem}
\newtheorem{lemma}{Lemma}

\newtheorem{algorithm}{Algorithm}
\newtheorem{proofprop}{\textit{Proof of Proposition}}
\newtheorem{remark}{Remark}
\newtheorem{prop}{Proposition}
\newtheorem{prooflemma}{\textit{Proof of Lemma}}

\newcommand\eqdef{\stackrel{\mathclap{\normalfont\mbox{def}}}{=}}
\newcommand*{\QED}{\hfill\ensuremath{\blacksquare}}
\DeclareMathOperator*{\argmin}{argmin} 
\newtheorem{propA}{Proposition}

\newtheorem{proofA}{\textit{Proof of Proposition}}

\newtheorem{proofthm}{Proof of Theorem}

\newtheorem{model}{\textbf{Model}}

\newcommand\tab[1][1cm]{\hspace*{#1}}

\begin{document}

\begin{center}
	\Large{{Gradient COBRA: A Kernel-based Consensual Aggregation for Regression}}\\
	
	\bigskip
	
	\normalsize
	Sothea Has
\end{center}

\begin{flushleft}
LPSM, Sorbonne Université Pierre et Marie Curie (Paris 6)\\ 75005 Paris, France\\
\url{sothea.has@lpsm.paris}
\end{flushleft}

\begin{abstract}
In this article, we introduce a kernel-based consensual aggregation method for regression problems. We aim to flexibly combine individual regression estimators $r_1, ..., r_M$ using a weighted average where the weights are defined based on predicted features given by all the basic estimators and some kernel function. This work extends the context of \cite{cobra} to a more general kernel-based framework. We show that this more general configuration also inherits the consistency of the basic consistent estimators, and the same convergence rate as in the classical method is achieved. Moreover, an optimization method based on gradient descent algorithm is proposed to efficiently and rapidly estimate the key parameter of the strategy. Various numerical experiments carried out on several simulated and real datasets are also provided to illustrate the efficiency and accuracy of the proposed method. Moreover, a domain adaptation-like property of the aggregation strategy is also illustrated on a physics data provided by Commissariat à l'Énergie Atomique (CEA).
\end{abstract}

\noindent \emph{Keywords:}
Consensual aggregation, kernel, regression.

\bigskip

\noindent \emph{2010 Mathematics Subject Classification:} {62G08, 62J99, 62P30}

\section{Introduction}
Aggregation methods, given the high diversity of available estimation strategies, are now of great interest in constructing predictive models. To this goal, several aggregation methods consisting of building a linear or convex combination of a collection of initial estimators have been introduced, for instance, in \cite{Cat04}, \cite{JN00}, \cite{Nem00}, \cite{yang2000,yang2001,yang2004}, \cite{bookDistributionFree}, \cite{W03}, \cite{Aud}, \cite{BTW06,BTW07a,BTW07b}, and \cite{DalTsy}. Other than aggregating, another possible approach is selecting the best estimator among the candidate estimators which is known as model selection technique (see, for example, \cite{MassStF}).

Apart from the usual linear combination and model selection methods, a different technique has been introduced in classification problems by \cite{mojirsheibani1999}. In his paper, the combination is the \textit{majority vote} among all the points for which their predicted classes, given by all the basic classifiers, \textit{coincide} with the predicted classes of the query point. Roughly speaking, instead of predicting a new point based on the structure of the original input, we look at the topology defined by the predictions of the candidate estimators. Each estimator was constructed differently so it may be able to capture different features of the input data and be useful in defining ``closeness''. Consequently, two points having similar predictions or classes seem reasonably having similar actual response values or belonging to the same actual class.

Later, \cite{mojirsheibani2000} and \cite{majidAndKong2016} introduced exponential and general kernel-based versions of the primal idea to improve the smoothness in selecting and weighting individual data points in the combination. In this context, the kernel function transforms the level of \textit{disagreements} between the predicted classes of a training point $x_i$ and the query point $x$ into a contributed weight given to the corresponding point in the vote. Besides, \cite{cobra} configured the original idea of \cite{mojirsheibani1999} as a regression framework where a training point $x_i$ is ``close'' to the query point $x$ if each of their predictions given by all the basic regression estimators is ``close''. Each of the close neighbors of $x$ will be given a uniformly 0-1 weight contributing to the combination. It was shown theoretically in these former papers that the combinations inherit the consistency property of consistent basic estimators.

Recently from a practical point of view, a kernel-based version of \cite{cobra} called \texttt{KernelCobra} has been implemented in \texttt{pycobra} python library (see \cite{pycobra}). This method has also been applied in filtering to improve the image denoising (see \cite{imageDenoisng}). Moreover, consensual aggregation methods such as \cite{cobra}, \cite{mixcobra} and the present method are also incorporated in a three-step methodology called \texttt{KFC procedure}, which combines unsupervised clustering and supervised prediction for (energy) data modeling (see \cite{kfc2021}). Such an idea of consensual aggregation was also used in unsupervised classification known as \texttt{Clustering Aggregation} (see, for example, \cite{clsuteraggregation2005} and \cite{clsuteraggregation2012}). On top of that, the aggregation method can also be used to handle the parameter tuning problem when different types of estimators are considered. It has been shown in \cite{has2022} that the method also maintains its good performance on highly correlated high-dimensional features of predictions that are plainly constructed without model selection or cross-validation.

In a complementary manner to the earlier works, we present in this paper a kernel-based consensual regression aggregation method, as well as its theoretical and numerical performances. More precisely, we show that the consistency inheritance property shown in \cite{cobra} also holds for this kernel-based configuration for a broad class of regular kernels. Moreover, evidence of numerical simulation carried out on several simulated models, and some real datasets, shows that the present method outperforms the classical one in both accuracy and efficiency.

This paper is organized as follows. Section~\ref{sec:theorem} introduces some notation, the definition of the proposed method, and presents the theoretical results, namely consistency and convergence rate of the variance-type term of the aggregation strategy. An optimization method based on gradient descent algorithm for estimating the bandwidth parameter is described in Section~\ref{sec:optimization}. Section~\ref{sec:numeric} illustrates the performances of the proposed method through several numerical experiments computed on different simulated and real datasets. Next, the conclusion and perspective, followed by the reproducibility of this study are given in Section~\ref{sec:Conclude} and Section~\ref{sec:supp} respectively. Lastly, Section~\ref{sec:proof} collects all the proofs of the theoretical results given in Section~\ref{sec:theorem}.

\section{The kernel-based combining regression}
\label{sec:theorem}
\subsection{Notation}
We consider a training sample $\mathcal{D}_n=\{(X_i,Y_i)_{i=1}^n\}$ where $(X_i,Y_i), i=1,2,...,n$, are {\it iid} copies of the generic couple $(X,Y)$. We assume that $(X,Y)$ is an $\mathbb{R}^d\times\mathbb{R}$-valued random variable with a suitable integrability which will be specified later. 

We randomly split the training data $\mathcal{D}_n$ into two parts of size $\ell$ and $k$ such that $\ell+k=n$. These are denoted by $\mathcal{D}_{\ell}=\{(X_i^{(\ell)},Y_i^{(\ell)})_{i=1}^{\ell}\}$ and $\mathcal{D}_k=\{(X_i^{(k)},Y_i^{(k)})_{i=1}^{k}\}$ respectively (a common choice is $k=\lceil n/2\rceil=n-\ell$). The $M$ basic regression estimators $r_{k,1},r_{k,2},...,r_{k,M}$ are constructed using only the data points in $\mathcal{D}_k$. These basic estimators can be any regression estimators such as linear regression, $k$NN, kernel smoother, SVR, lasso, ridge, neural networks, naive Bayes, bagging, gradient boosting, random forests, etc. They could be parametric, nonparametric or semi-parametric with their possible tuning parameters. For the combination, we only need the predictions given by all these basic estimators of the remaining part $\mathcal{D}_{\ell}$ and  the query point $x$. 

In the sequel, for any $x\in\mathbb{R}^d$, the following notation is used:
\begin{itemize}
	\item$\textbf{r}_k(x)=(r_{k,1}(x),r_{k,2}(x),...,r_{k,M}(x))$: the vector of predictions of $x$.
	\item$\|x\|=\|x\|_2=\sqrt{\sum_{i=1}^dx_i^2}$: Euclidean norm on $\mathbb{R}^d$.
	\item$\|x\|_1=\sum_{i=1}^d|x_i|$: $\ell_1$ norm on $\mathbb{R}^d$.
	\item$g^*(x)=\mathbb{E}[Y|X=x]$: the regression function.
	\item$g^*(\textbf{r}_k(x))=\mathbb{E}[Y|\textbf{r}_k(x)]$: the conditional expectation of the response variable given all the predictions. This can be proven to be the optimal estimator in regression over the set of predictions $\textbf{r}_k(X)$.
	\item $\mathds{1}_{\{p\}}=\begin{cases}1,&\mbox{if }p\mbox{ is true}\\0,&\mbox{otherwise}\end{cases}$: the indicator function.
\end{itemize}
The consensual regression aggregation is the weighted average defined by
\begin{equation}
	\label{eq:combine}
	g_n(\textbf{r}_k(x))=\sum_{i=1}^{\ell}W_{n,i}(x)Y_i^{(\ell)}.
\end{equation}

Recall that given all the basic estimators $r_{k,1},r_{k,2},...,r_{k,M}$, the aggregation method proposed by \cite{cobra} corresponds to the following naive weights:

\begin{equation}
	\label{eq:cobra1}
	W_{n,i}(x)=\frac{\prod_{m=1}^M\displaystyle\mathds{1}_{\{|r_{k,m}(X_i)-r_{k,m}(x)|<h\}}}{\sum_{j=1}^{\ell}\prod_{m=1}^M\mathds{1}_{\{|r_{k,m}(X_j)-r_{k,m}(x)|<h\}}},i=1,2,...,\ell.
\end{equation}
Moreover, the condition of ``closeness for all'' predictions, can be relaxed to ``some'' predictions, which corresponds to the following weights:
\begin{equation}
	\label{eq:cobra2}
	W_{n,i}(x)=\frac{\displaystyle\mathds{1}_{\{\sum_{m=1}^M\mathds{1}_{\{|r_{k,m}(X_i)-r_{k,m}(x)|<h\}}\geq \alpha M\}}}{\sum_{j=1}^{\ell}\mathds{1}_{\{\sum_{m=1}^M\mathds{1}_{\{|r_{k,m}(X_j)-r_{k,m}(x)|<h\}}\geq \alpha M\}}},i=1,2,...,\ell
\end{equation}
where $\alpha\in\{1/M,2/M,...,1\}$ is the proportion of consensual predictions required and $h>0$ is the bandwidth or window parameter to be determined. Constructing the proposed method is equivalent to searching for the best possible value of these parameters over a given grid, minimizing some quadratic error which will be described in Section~\ref{sec:optimization}.

In the present paper, $K:\mathbb{R}^M\to\mathbb{R}_+$ denotes a regular kernel which is a decreasing function satisfying:
\begin{align}
	\label{eq:regular}
	\exists b,\kappa_0,\rho>0\ \text{such that}
	\begin{cases}
		b\mathds{1}_{B_M(0,\rho)}(z)\leq K(z)\leq 1, \forall z\in\mathbb{R}^M\\
		\int_{\mathbb{R}^M}\sup_{u\in B_M(z,\rho)}K(u)dz = \kappa_0 < +\infty
	\end{cases}
\end{align}
where $B_M(c,r)=\{z\in\mathbb{R}^M:\|c-z\|<r\}$ denotes the open ball of center $c\in\mathbb{R}^M$ and radius $r>0$ of $\mathbb{R}^M$. We propose in equation~\eqref{eq:combine} a method associated to the weights defined at any query point $x\in\mathbb{R}^d$ by
\begin{align}
	\label{eq:KCOBRA}
	W_{n,i}(x)=\frac{K_h(\textbf{r}_k(X_i^{(\ell)})-\textbf{r}_k(x))}{\sum_{j=1}^{\ell}K_h(\textbf{r}_k(X_j^{(\ell)})-\textbf{r}_k(x))},i=1,2,...,\ell 
\end{align}
where $K_h(z)=K(z/h)$ for some bandwidth parameter $h>0$ with the convention of $0/0=0$. Observe that the combination in equation~\eqref{eq:combine} is computed based only on $\mathcal{D}_{\ell}$ but the construction of the method depends on the whole training data $\mathcal{D}_n$ as the basic estimators are all constructed using $\mathcal{D}_k$. In our setting, we treat the vector of predictions $\textbf{r}_k(x)$ as an $M$-dimensional feature, and the kernel function is applied on the whole vector at once. Note that the implementation of \texttt{KernelCobra} in \cite{Guedj_2020} corresponds to the following weights:
\begin{equation}
	\label{eq:Guedj}
	W_{n,i}(x)=\frac{\sum_{m=1}^MK_h(r_{k,m}(X_i^{(\ell)})-r_{k,m}(x))}{\sum_{j=1}^{\ell}\sum_{m=1}^MK_h(r_{k,m}(X_j^{(\ell)})-r_{k,m}(x))},i=1,2,...,\ell
\end{equation}
where the univariate kernel function $K$ is applied on each component of the predicted vector $\textbf{r}_k(.)$ separately. 
In this case, the weight $W_{n,i}(x)$ defined in equation~\eqref{eq:Guedj} above is more costly in computing than the one in the proposed method since the univariate kernel function has to be applied on all the entries of vectors $\textbf{r}_k(X_i^{\ell})-\textbf{r}_k(x)=(r_{k,1}(X_i^{\ell})-r_{k,1}(x),...,r_{k,M}(X_i^{\ell})-r_{k,M}(x))$ for all $i=1,...,\ell$. This entry-wise operation prevents us from trading off memory storage for computational complexity. On the other hand, the weights in equation~\eqref{eq:KCOBRA} of the proposed method depend on pair-wise distances between the predicted vectors of the training points $X_i^{(\ell)}$'s and the query point $x$, $d'(\textbf{r}_k(X_i),\textbf{r}_k(x))$, for some distance $d'$ (associated to the kernel function). This dependency allows us to trade the memory storage off for computational complexity, yielding more efficient computation and the implementation of an optimization procedure based on gradient descent algorithm (section~\ref{sec:optimization}). 

\subsection{Theoretical performance}
The performance of the combining estimation $g_n$ is measured using the quadratic risk defined by
\begin{align*}
	\mathbb{E}\Big[|g_n(\textbf{r}_k(X))-g^*(X)|^2\Big]
\end{align*}
where the expectation is taken with respect to both $X$ and the training sample $\mathcal{D}_n$. Firstly, we begin with a simple decomposition of the distortion between the proposed method and the optimal regression estimator $g^*(X)$ by introducing the optimal regression estimator over the set of predictions $g^*(\textbf{r}_k(X))$. The following proposition shows that the nonasymptotic-type control of the distortion, presented in Proposition.2.1 of \cite{cobra}, also holds for this case of regular kernels.

\begin{prop}
	\label{prop:1}
	Let $\textbf{r}_k=(r_{k,1},r_{k,2},...,r_{k,M})$ be the collection of all basic estimators, and let $g_n(\textbf{r}_k(x))$ be the combined estimator defined in equation~\eqref{eq:combine} with the weights given in equation~\eqref{eq:KCOBRA} computed at point $x\in\mathbb{R}^d$. Then, for all distributions of $(X,Y)$ with $\mathbb{E}[|Y|^2]< +\infty$,
	\begin{align*}
		\mathbb{E}\Big[|g_n(\textbf{r}_k(X))-g^*(X)|^2\Big]&\leq \inf_{f\in\mathcal{G}}\mathbb{E}\Big[|f(\textbf{r}_k(X))-g^*(X)|^2\Big]\\
		&\quad+\mathbb{E}\Big[|g_n(\textbf{r}_k(X))-g^*(\textbf{r}_k(X))|^2\Big]
	\end{align*}
	where $\mathcal{G}$ is the class of any function $f:\mathbb{R}^M\to\mathbb{R}$ satisfying $\mathbb{E}[f(\textbf{r}_k(X))|^2]<+\infty$. In particular,
	\begin{align*}
		\mathbb{E}\Big[|g_n(\textbf{r}_k(X))-g^*(X)|^2\Big]&\leq \min_{1\leq m\leq M}\mathbb{E}\Big[|r_{k,m}(X)-g^*(X)|^2\Big]\\
		&\quad+\mathbb{E}\Big[|g_n(\textbf{r}_k(X))-g^*(\textbf{r}_k(X))|^2\Big].
	\end{align*}
\end{prop}


\noindent The two terms of the last bound can be viewed as a bias-variance decomposition where the first term $\min_{1\leq m\leq M}\mathbb{E}[|r_{k,m}(X)-g^*(X)|^2]$ can be seen as the bias and $\mathbb{E}[|g_n(\textbf{r}_k(X))-g^*(\textbf{r}_k(X))|^2]$ is the variance-type term (\cite{cobra}). Given all the estimators, the first term cannot be controlled as it depends on the performance of the best constructed estimator, and it will be the asymptotic performance of the proposed method. Our main task is to deal with the second term, which can be proven to be asymptotically negligible in the following key proposition.

\begin{prop}
	\label{prop:2}
	Assume that $r_{k,m}$ is bounded for all $m=1,2,..,M$. Let $h\rightarrow0$ and $\ell\rightarrow+\infty$ such that $h^M\ell\to+\infty$. Then
	\begin{align*}
		\mathbb{E}\Big[|g_n(\textbf{r}_k(X))-g^*(\textbf{r}_k(X))|^2\Big]\rightarrow0\ \text{as }\ell\rightarrow+\infty
	\end{align*}
	for all distribution of $(X,Y)$ with $\mathbb{E}[|Y|^2]<+\infty$. Thus,
	\begin{align*}
		\limsup_{\ell\rightarrow+\infty}\mathbb{E}\Big[|g_n(\textbf{r}_k(X))-g^*(X)|^2\Big]\leq\inf_{f\in\mathcal{G}}\mathbb{E}\Big[|f(\textbf{r}_k(X))-g^*(X)|^2\Big].
	\end{align*}
	And in particular,
	\begin{align*}
		\limsup_{\ell\rightarrow+\infty}\mathbb{E}\Big[|g_n(\textbf{r}_k(X))-g^*(X)|^2\Big]\leq\min_{1\leq m \leq M}\mathbb{E}\Big[|r_{k,m}(X)-g^*(X)|^2\Big].
	\end{align*}
\end{prop}
Proposition~\ref{prop:2} above is an analogous setup of Proposition 2.2 in \cite{cobra}. To prove this result, we follow the procedure of Stone's theorem (see, for example, \cite{stone1977} and Chapter 4 of \cite{bookDistributionFree}) of weak universal consistency of non-parametric regression. However, showing this result for the class of regular kernels is not straightforward. Most of the previous studies provided such a result of $L_2$-consistency only for the class of compactly supported kernels (see, for example, Chapter 5 of \cite{bookDistributionFree}). In this study, we can derive the result for this broader class thanks to the boundedness of all basic estimators. However, the price to pay for the universality for this class of regular kernels is the lack of convergence rate. To this goal, a weak smoothness assumption of $g^*$ with respect to the basic estimators is required. For example, the convergence rate of the variance-type term in \cite{cobra} is of order $O(\ell^{-2/(M+2)})$ under the same smoothness assumption, and this result also holds for all the compactly support kernels. In this study, we can derive the same convergence rate for the class of kernel functions with the tails increase at least of exponential speed. This main theoretical result is given in the following theorem.

\begin{theorem}
	\label{thm.1}
	Assume that the response variable $Y$ and all the basic estimators $r_{k,m},m=1,2,...,M$, are bounded by some constant $R$. Suppose that there exists a constant $L\geq0$ such that, for every $k\geq1$,
	\begin{equation*}
		|g^*(\textbf{r}_k(x))-g^*(\textbf{r}_k(y))|\leq L\|\textbf{r}_k(x)-\textbf{r}_k(y)\|,\forall x,y\in\mathbb{R}^d.
	\end{equation*}
	We assume moreover that there exists some positive constants $\alpha,R_K$ and $C_K$ such that
	\begin{equation}
		\label{eq:assumption}
		K(z)\leq C_K\exp(-\|z\|^{\alpha}), \forall z\in\mathbb{R}^M, \|z\|\geq R_K.
	\end{equation}
	Then,  one has
	\begin{equation}
		\label{eq:convergRate}
		\mathbb{E}[|g_n(\textbf{r}_k(X))-g^*(X)|^2]\leq \min_{1\leq m\leq M}\mathbb{E}[|r_{k,m}(X)-g^*(X)|^2]+C\ell^{-\frac{2}{M+2}}
	\end{equation}
	for some positive constant $C=C(b,L,R,R_K,C_K)$ independent of $\ell$. 
\end{theorem} 

\noindent From this result, if there exists a consistent estimator named $r_{k,m_0}$ in the list $\{r_{k,m}\}_{m=1}^M$ i.e.,
$$\mathbb{E}[|r_{k,m_0}(X)-g^*(X)|^2]\to0\ \ \text{as }k\to+\infty,$$
then the combining estimator $g_n$ is also consistent for all distribution of ($X,Y$) in some class $\mathcal{M}$. Consequently, under the assumption of Theorem~\ref{thm.1}, one has 
$$\lim_{k,\ell\to+\infty}\mathbb{E}[|g_n(\textbf{r}_k(X))-g^*(X)|^2]=0.$$

\section{Bandwidth estimation using gradient descent}
\label{sec:optimization}

In earlier works by \cite{cobra} and \cite{Guedj_2020}, the training data $\mathcal{D}_{n}$ is practically broken down into three balanced parts: $\mathcal{D}_k$ for constructing all candidate estimators $\{\textbf{r}_{k,m}\}_{m=1}^M$, $\mathcal{D}_{\ell_1}$ for building aggregation defined in equation~\eqref{eq:combine}, and $\mathcal{D}_{\ell_2}$ for tuning the key parameters of the methods. Within these frameworks, the bandwidth parameter $h$ is estimated by minimizing the following loss,
\begin{align}
	\label{eq:error1}
	\varphi_M(h)=\frac{1}{|\mathcal{D}_{\ell_2}|}\sum_{(X_j,Y_j)\in\mathcal{D}_{\ell_2}}[g_n(\textbf{r}_k(X_j))-Y_j]^2,
\end{align}
where $|\mathcal{D}_{\ell_2}|$ denotes the cardinality of $\mathcal{D}_{\ell_2}$, and $g_n(\textbf{r}_k(X_j))=\sum_{(X_i,Y_i)\in\mathcal{D}_{\ell_1}}W_{n,i}(X_j)Y_i$ is given in equation~\eqref{eq:combine}. Note that the subscript $M$ of $\varphi_M(h)$ indicates the full consensus between the $M$ components of the predictions $\textbf{r}_k(X_i)$ and $\textbf{r}_k(X_j)$ for any $X_i$ of $\mathcal{D}_{\ell_1}$ and $X_j$ of $\mathcal{D}_{\ell_2}$. In this case, constructing an aggregation method $g_n$ is equivalent to searching for an optimal parameter $h^*$ over a given grid $\mathcal{H}=\{h_{\min},...,h_{\max}\}$ i.e.,
$$h^*=\argmin_{h\in\mathcal{H}}\varphi_M(h).$$ 
The parameter $\alpha$ of equation~\eqref{eq:cobra2} can be tuned easily by considering $\varphi_{\alpha M}(h)$ where $\alpha\in\{1/M,2/M,...,1\}$ referring to the proportion of consensuses required among the $M$ components of the predictions. In this case, the optimal parameters $\alpha^*$ and $h^*$ are chosen to be the minimizer of $\varphi_{\alpha M}(h)$ i.e.,
$$(\alpha^*,h^*)=\argmin_{(\alpha,h)\in\{1/M,2/M,...,1\}\times\mathcal{H}}\varphi_{\alpha M}(h).$$ 
Note that in both papers, the grid search algorithm is used in searching for the optimal bandwidth parameter. 

In this paper, the training data is broken down into only two parts, $\mathcal{D}_{k}$ and $\mathcal{D}_{\ell}$. Again, we construct the basic estimators using $\mathcal{D}_{k}$, and for any $\kappa$ folds $F_1,...,F_{\kappa}$ ($\kappa\geq 2$) of $\mathcal{D}_{\ell}$, we propose the following $\kappa$-fold cross-validation error which is a function of the bandwidth parameter $h>0$ defined by
\begin{align}
	\label{eq:error2}
	\varphi^{\kappa}(h)=\frac{1}{\kappa}\sum_{p=1}^{\kappa}\sum_{(X_j,Y_j)\in F_p}[g_n(\textbf{r}_k(X_j))-Y_j]^2
\end{align}
where in this case, $g_n(\textbf{r}_k(X_j))=\sum_{(X_i,Y_i)\in \mathcal{D}_{\ell}\setminus F_p}W_{n,i}(X_j)Y_i$, is computed using the remaining $\kappa-1$ folds of $\mathcal{D}_{\ell}$ leaving $F_p\subset\mathcal{D}_{\ell}$ as the corresponding validation fold\footnote{In this part, we simply write $(X_i,Y_i)\in\mathcal{D}_{\ell}$ without the superscript $(\ell)$.}. We often observe the convex-like curves of the cross-validation quadratic error on many simulations; and from this observation, we propose using a gradient descent algorithm to estimate the optimal bandwidth parameter. The associated gradient descent algorithm used to estimate the optimal parameter $h^*$ is implemented as follows:
\begin{tcolorbox}
	\begin{algorithm}{: Gradient descent for estimating $h^*$:}
		\begin{enumerate}
			\item Initialization: $h_0$, a learning rate $\lambda>0$, threshold $\delta>0$ and the maximum number of iteration $N$.
			\item For $k=1,2,...,N$, \textbf{while} $\Big|\frac{d}{d h}\varphi^{\kappa}(h_{k-1})\Big|>\delta$ do: 
			$$h_k\gets h_{k-1}-\lambda\frac{d}{d h}\varphi^{\kappa}(h_{k-1})$$
			\item return $h_k$ violating the \textbf{while} condition or $h_N$ to be the estimation of $h^*$.
		\end{enumerate}
	\end{algorithm}
\end{tcolorbox}
\noindent From equation~\eqref{eq:error2}, for any $(X_j,Y_j)\in F_p$, one has
\small{
	\begin{align*}
		\frac{d}{d h}\varphi^{\kappa}(h)&=\frac{1}{\kappa}\sum_{p=1}^{\kappa}\sum_{(X_j,Y_j)\in F_p}2\frac{\partial}{\partial h}g_n(\textbf{r}_k(X_j))(g_n(\textbf{r}_k(X_j))-Y_j)
	\end{align*}
}
where 
\small{
	\begin{align*}
		g_n(\textbf{r}_k(X_j))&=\frac{\sum_{(X_i,Y_i)\mathcal{D}_{\ell}\in\setminus F_p}Y_iK_h(\textbf{r}_k(X_j)-\textbf{r}_k(X_i))}{\sum_{(X_q,Y_q)\in\mathcal{D}_{\ell}\setminus F_p}K_h(\textbf{r}_k(X_j)-\textbf{r}_k(X_q))}.
	\end{align*}
	This implies that
	\begin{align*}
		\frac{\partial}{\partial h}g_n(\textbf{r}_k(X_j))&=\sum_{(X_i,Y_i),(X_q,Y_q)\in\mathcal{D}_{\ell}\setminus F_p}(Y_i-Y_q)\frac{\frac{\partial}{\partial h}K_h(\textbf{r}_k(X_j)-\textbf{r}_k(X_i))K_h(\textbf{r}_k(X_j)-\textbf{r}_k(X_q))}{\Big[\sum_{(X_i,Y_i)\mathcal{D}_{\ell}\in\setminus F_p}K_h(\textbf{r}_k(X_j)-\textbf{r}_k(X_i))\Big]^2}.
	\end{align*}
}

\noindent The differentiability of $g_n$ depends entirely on the kernel function $K$. Therefore, for suitable kernels, the implementation of the algorithm is straightforward. For example, in the case of Gaussian kernel $K_h(x)=\exp(-h\|x\|^2/(2\sigma^2))$ for some $\sigma>0$, one has
\small{
	\begin{align*}
		\frac{\partial}{\partial h}g_n(\textbf{r}_k(X_j))&=\sum_{(X_i,Y_i),(X_q,Y_q)\in\mathcal{D}_{\ell}\setminus F_p}(Y_q-Y_i)\|\textbf{r}_k(X_j)-\textbf{r}_k(X_i)\|^2\times\\
		&\quad\frac{\exp\Big(-h(\|\textbf{r}_k(X_j)-\textbf{r}_k(X_i)\|^2+\|\textbf{r}_k(X_j)-\textbf{r}_k(X_q)\|^2)/(2\sigma^2)\Big)}{2\sigma^2\Big(\sum_{(X_q,Y_q)\notin F_p}\exp(-h\|\textbf{r}_k(X_j)-\textbf{r}_k(X_q)\|^2/(2\sigma^2))\Big)^2}.
	\end{align*}
}
This suggests that we only need to store the distance matrices $D_{p}=(d'_{qj})$ where $d'_{qj}=\|\textbf{r}_k(X_q)-\textbf{r}_k(X_j)\|^2$ is the squared Euclidean distance between predictions of the input data from the $\kappa -1$ folds $\mathcal{D}_{\ell}\setminus F_p$ and the corresponding validation fold $F_p$ for $p=1,...,\kappa$. Then, the gradient can be computed straight away for any smoothing parameter $h>0$. 

To prevent the algorithm from reaching negative values of the bandwidth parameter during operation, a few adjustments have been implemented. Firstly, the predicted features are normalized for example, to be in the range $[0,1]^M$. Then, the error is computed at a few randomly selected bandwidth parameters, and the algorithm begins at the parameter with the lowest error. Additionally, the learning rate $\lambda$ is decreased when the parameter takes negative values, which may occur due to a large learning rate. To handle cases where the error curve is very flat around the optimal bandwidth, an option has been included to adjust the speed of the learning rate. This approach has resulted in faster algorithm performance, without requiring knowledge of the interval containing the optimal parameter, as with grid search. Moreover, it is possible to estimate the parameter that causes the gradient of the objective function to vanish. This leads to a well-constructed aggregation method, as reported in the next section.

\section{Numerical examples}
\label{sec:numeric}
This section is devoted to numerical experiments to illustrate the performance of our proposed method. It is shown in \cite{cobra} that the classical method mostly outperforms the basic estimators of the combination. In this experiment, we compare the performances of the proposed methods with the classical one and all the basic regressors. Several options of kernel functions are considered. Most kernels are compactly supported on $[-1,1]$, taking nonzero values only on $[-1,1]$, except for the case of compactly supported Gaussian which is supported on $[-\rho_1,\rho_1]$, for some $\rho_1>0$. Moreover to implement the gradient descent algorithm in estimating the bandwidth parameter, we also present the results of non-compactly supported cases such as classical Gaussian and 4-exponential kernels. All kernels considered in this paper are listed in Table~\ref{tab:kernels}, and some of them are displayed (univariate case) in Figure~\ref{fig:1} below.

\begin{table}[h!]
	\centering 
	\small
	\begin{tabular}{| l | l |}  
		\hline                      
		\tab[0.45cm] \textbf{Kernel} & \tab[1.3cm] \textbf{Formula}\\
		\hline
		Naive\footnote{The naive kernel corresponds to the method by \cite{cobra}.} & $K(x)=\prod_{i=1}^d\mathds{1}_{\{|x_i|\leq 1\}}$\\ 
		\hline
		Epanechnikov & $K(x)=(1-\|x\|^2)\mathds{1}_{\{\|x\|\leq 1\}}$\\
		\hline 
		Bi-weight & $K(x)=(1-\|x\|^2)^2\mathds{1}_{\{\|x\|\leq 1\}}$\\
		\hline
		Tri-weight & $K(x)=(1-\|x\|^2)^3\mathds{1}_{\{\|x\|\leq 1\}} $\\
		\hline
		Compact-support Gaussian & $K(x)=\exp\{-\|x\|^2/(2\sigma^2)\}\mathds{1}_{\{\|x\|\leq \rho_1\}}, \sigma,\rho_1>0$\\
		\hline                          
		Gaussian & $K(x)=\exp\{-\|x\|^2/(2\sigma^2)\}, \sigma>0$\\
		\hline  
		$4$-exponential & $K(x)=\exp\{-\|x\|^4/(2\sigma^4)\}, \sigma>0$\\
		\hline
	\end{tabular}
	\caption{Kernel functions used.}%
	\label{tab:kernels}
\end{table}

\begin{figure}[h!]
	\centering
	\begin{tikzpicture}
		\begin{axis}[
			axis lines = left,
			xlabel = $x$,
			ylabel = {$K(x)$},
			legend style={at={(1,1)},anchor=north west}, 
			legend cell align={left},
			width=10cm,height=6.5cm
			]
			\addplot [
			domain=-3.5:-1, 
			samples=100, 
			color=black,
			line width=1.5pt,
			]
			{0};
			\addplot [
			domain=1:3.5, 
			samples=100, 
			color=black,
			line width=1.5pt,
			]
			{0};
			
			\addplot [
			domain=-1:1, 
			samples=100, 
			color=black,
			line width=1.5pt,
			]
			{1};
			\addplot [
			densely dashed,
			domain=-1:1, 
			samples=100, 
			color=blue,
			line width=1pt
			]
			{1-x^2};
			\addplot [
			domain=-1:1, 
			samples=100, 
			color=orange,
			line width=1.1pt
			]
			{(1-x^2)^2};
			\addplot [
			densely dashdotted,
			domain=-1:1, 
			samples=100, 
			color=purple,
			line width=1.2pt
			]
			{(1-x^2)^3};
			\addplot [
			loosely dashdotted,
			domain=-3.5:3.5, 
			samples=100, 
			color=red,
			line width=1.4pt
			]
			{exp(-x^2)};
			\addplot [
			loosely dotted,
			domain=-3.5:3.5, 
			samples=100, 
			color=cyan,
			line width=1.5pt
			]
			{exp(-x^4)};
			\legend{,,Naive,Epanechnikov,Bi-weight,Tri-weight,Gaussian, 4-exponential}
		\end{axis}
	\end{tikzpicture}
	\caption{The shapes of some kernels.}
	\label{fig:1}
\end{figure}
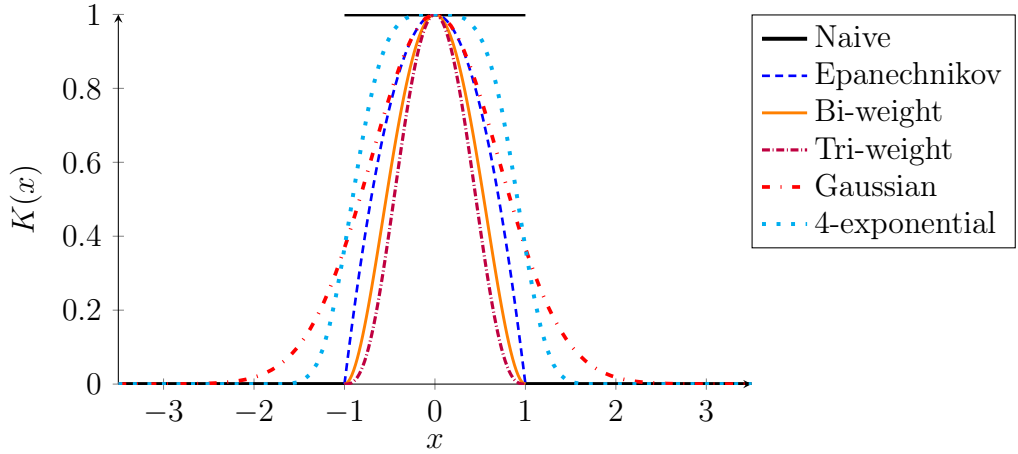

\subsection{Simulated datasets}
\label{subsec:simulated}
In this subsection, we study the performances of our proposed method on the same set of simulated datasets of size $n$ as provided in \cite{cobra}. The input data is either independent and uniformly distributed over $(-1,1)^d$ ({\it uncorrelated} case) or distributed from a Gaussian distribution $\mathcal{N}(0,\Sigma)$ where the covariance matrix $\Sigma$ is defined by $\Sigma_{ij}=2^{-|i-j|}$ for $1\leq i,j\leq d$ ({\it correlated} case). We consider the following models:

\begin{model}{\bf:}
	\label{mod:1}
	$n=800, d=50, Y=X_1^2+\exp(-X_2^2).$
\end{model}
\begin{model}{\bf:}
	\label{mod:2}
	$n=600, d=100, Y=X_1X_2+X_3^2-X_4X_7+X_8X_{10}-X_6^2+\mathcal{N}(0,0.5).$
\end{model}
\begin{model}{\bf:}
	\label{mod:3}
	$n=600, d=100, Y=-\sin(2X_1)+X_2^2+X_3-\exp(-X_4)+\mathcal{N}(0,0.5).$
\end{model}
\begin{model}{\bf:}
	\label{mod:4}
	$n=600, d=100, Y=X_1+(2X_2-1)^2+\sin(2\pi X_3)/(2-\sin(2\pi X_3))+\sin(2\pi X_4)+2 \cos(2\pi X_4)+3\sin^2(2\pi X_4)+4\cos^2(2\pi X_4)+\mathcal{N}(0,0.5).$
\end{model}
\begin{model}{\bf:}
	\label{mod:5}
	$n=700, d=20, Y=\mathds{1}_{\{X_1>0\}}+X_2^3+\mathds{1}_{\{X_4+X_6-X_8-X_9>1+X_{14}\}}+\exp(-X_2^2)+\mathcal{N}(0,0.05).$
\end{model}
These first five models are taken from \cite{cobra} which allows us to compare the performance of the methods. Note that by the design, there are not many active predictors contributing to the target, and most of them act as the noise. To see how the proposed method behaves on different type of datasets where more active independent variables are presented, we introduce the following models:
\begin{model}{\bf:}
	\label{mod:6}
	$n=500, d=20, Y=(\sum_{j=1}^5\sum_{k=0}^3X_{j+5k})\cos((\prod_{k=1}^5X_{4k})\pi/2)+\mathcal{N}(0, 0.25)$
\end{model}
\begin{model}{\bf:}
	\label{mod:7}
	$n=600, d=30, Y= \sum_{j=1}^{15}e^{0.25-X_j^2}\sin(\pi X_{j+15})+\mathcal{N}(0,0.25)$
\end{model}
\begin{model}{\bf:}
	\label{mod:8}
	$n=700, d=50, Y=(\sum{j=1}^25X_{2j}\sin(\pi/X_{2j-1}))e^{\sum{k=1}^5X_{10k}^2/10}+\mathcal{N}(0,0.75)$
\end{model}
Moreover, it is interesting to consider some high-dimensional cases as many real problems such as image and signal processing involve these kinds of datasets. Therefore, we also consider the following two high-dimensional models where all the independent variables contribute to the target via the coefficient $\beta_j$'s.
\begin{model}{\bf:}
	\label{mod:9}
	$n=600, d=1500, Y=\pi+\sum_{j=1}^d\beta_j\frac{X_j\log|5+X_j|}{1+e^{X_j}}+\mathcal{N}(0,1)$,\ \textrm{where }$\beta_j=2^{-(d+1-j)/50}+3^{-j/50}, j = 1,...,d$.
\end{model}
\begin{model}{\bf:}
	\label{mod:10}
	$n=700, d=1500, Y=e+\sum_{j=1}^d\beta_j\frac{X_je^{-X_j}}{1-\log|10-X_j|}+\mathcal{N}(0,1.25)$,\ \textrm{where }$\beta_j=\frac{e^{-j/30}}{1-e^{-(d+1-j)/30}}, j = 1,...,d$.
\end{model}
For each model, the proposed method is implemented over $100$ replications. We randomly split $80\%$ of each simulated dataset into two equal parts, $\mathcal{D}_{\ell}$ and $\mathcal{D}_k$ where $\ell=\lceil 0.8\times n/2\rceil-k$, and the remaining $20\%$ is treated as the corresponding testing data. We measure the performance of any regression method $f$ using {\it root mean square error} (RMSE) evaluated on the $20\%$-testing data defined by

\begin{equation}
	\text{RMSE}(f)=\left(\frac{1}{n_{\text{test}}}\sum_{i=1}^{n_{\text{test}}}(y_i^{\text{test}}-f(x_i^{\text{test}}))^2\right)^{1/2}.
\end{equation}

Table~\ref{tab:uncorr} and \ref{tab:corr} below contain the average RMSEs and the corresponding standard errors (into brackets) over $100$ runs of {\it uncorrelated} and {\it correlated} cases respectively. In each table, the first block contains five columns corresponding to the following five basic regressors ${\bf r}_k=(r_{k,m})_{m=1}^5$:
\begin{itemize}
	\item \textcolor{cyan}{\bf Rid}: Ridge regression (R package \texttt{glmnet}, see \cite{glmnet}).
	\item \textcolor{cyan}{\bf Las}: Lasso regression (R package \texttt{glmnet}).
	\item \textcolor{cyan}{\bf$k$NN}: $k$-nearest neighbors regression (R package \texttt{FNN}, see \cite{FNN}).
	\item \textcolor{cyan}{\bf Tr}: Regression tree (R package \texttt{tree}, see \cite{tree}).
	\item \textcolor{cyan}{\bf RF}: Random Forest regression (R package \texttt{randomForest}, see \cite{randomForest}).
\end{itemize}
We choose $k=5$ for $k$-NN and $ntree = 300$ for random forest algorithm, and other methods are implemented using the default parameters. The best performance of each method in this block is given in \textbf{boldface}. The second block contains the last eight columns corresponding to kernel functions and different types of aggregation methods. The abbreviations of all the methods in this block are given below:
\begin{itemize}
	\item \textcolor{cyan}{\bf COBRA}: the \textit{classical COBRA} by \cite{cobra}.
	\item \textcolor{cyan}{\bf Epan}: the aggregation method using \textit{Epanechnikov} kernel.
	\item \textcolor{cyan}{\bf Bi-wgt}: the aggregation method using \textit{Bi-weight} kernel.
	\item \textcolor{cyan}{\bf Tri-wgt}: the aggregation method using \textit{Tri-weight} kernel.
	\item \textcolor{cyan}{\bf C-Gaus}: the aggregation method using \textit{Compact Gaussian} kernel.
	\item \textcolor{cyan}{\bf Gauss}: the aggregation method using \textit{Gaussian} kernel. 
	\item \textcolor{cyan}{\bf Exp$4$}: the aggregation method using \textit{$4$-Exponential} kernel.
	\item \textcolor{cyan}{\bf KCOBRA}: the \textit{KernelCobra} by \cite{pycobra}.
\end{itemize}
The optimal RMSEs of each model in this block is also written in \textbf{boldface}. For all the compactly supported kernels, we consider $500$ values of bandwidth parameter $h$ in a uniform grid $\{10^{-100},...,h_{\max}\}$ where $h_{\max}=10$, which is chosen to be large enough, likely to contain the optimal parameter to be searched. For the compactly supported Gaussian kernel, we set $\rho_1=3$ and $\sigma=1$ therefore its support is $[-3,3]$. For the two non-compactly supported kernels, Gaussian and 4-exponential, the optimal parameters are estimated using gradient descent algorithm described in the previous section. Lastly, Gaussian kernel is used for \texttt{KernelCobra}, and the optimal bandwidth is estimated using \texttt{optimal\_kernelbandwidth} method of \texttt{pycobra} library.

\begin{sidewaystable}[ph!]
	\renewcommand{\arraystretch}{1.25}
	\tiny
	\centering  
	\caption{Average MSEs in the uncorrelated case.}
	\label{tab:uncorr}
	
	\hspace{0.5em}
	
	\begin{tabular}{ c | c c c c c | c c c c c c c c}  
		\hline                      
		\textcolor{cyan}{\bf Model} & \textcolor{cyan}{\bf Las} &\textcolor{cyan}{\bf Rid} &\textcolor{cyan}{\bf $k$NN} & \textcolor{cyan}{\bf Tr} & \textcolor{cyan}{\bf RF} & \textcolor{cyan}{\bf COBRA} & \textcolor{cyan}{\bf Epan} & \textcolor{cyan}{\bf Bi-wgt} & \textcolor{cyan}{\bf Tri-wgt}  & \textcolor{cyan}{\bf C-Gaus} & \textcolor{cyan}{\bf Gauss} & \textcolor{cyan}{\bf Exp$4$} & \textcolor{cyan}{\bf KCOBRA}\\ [0.5ex]   
		\hline  
		\multirow{2}{*}{\ref{mod:1}} & \multirow{2}{*}{\raisebox{2ex}{$0.156$}} & \multirow{2}{*}{\raisebox{2ex}{$0.133$}} & \multirow{2}{*}{\raisebox{2ex}{$0.143$}} & \multirow{2}{*}{\raisebox{2ex}{$\bf 0.027$}} & $\multirow{2}{*}{\raisebox{2ex}{0.032}}$ & $0.020$ & $0.018$ & $0.017$ & $0.017$ & $0.017$ & $\bf 0.015$ & $0.016$ & $0.061$\\
		& (0.016) & (0.013)& (0.014) & (0.004) & (0.004) & $(0.004)$ & $(0.003)$ & $(0.003)$ & $(0.003)$ & $(0.003)$ & $(0.002)$ & (0.003) & (0.027)\\
		\hline 
		\multirow{2}{*}{\ref{mod:2}} & \multirow{2}{*}{\raisebox{2ex}{$1.301$}} & \multirow{2}{*}{\raisebox{2ex}{$0.784$}} & \multirow{2}{*}{\raisebox{2ex}{$0.873$}} & \multirow{2}{*}{\raisebox{2ex}{$1.124$}} & $\multirow{2}{*}{\raisebox{2ex}{\textbf{0.707}}}$ & $0.722$ & $0.718$ & $0.712$ & $0.715$ & $0.712$ & $\textbf{0.709}$ & $0.710$ & $0.788$\\
		&(0.216) & (0.110)& (0.123) & (0.165) & (0.097) &$(0.065)$ & $(0.079)$ & $(0.080)$ & $(0.079)$ & $(0.079)$ & $(0.078)$ & (0.079) & (0.085)\\
		\hline  
		\multirow{2}{*}{\ref{mod:3}} & \multirow{2}{*}{\raisebox{2ex}{$0.664$}} & \multirow{2}{*}{\raisebox{2ex}{$0.669$}} & \multirow{2}{*}{\raisebox{2ex}{$1.477$}} & \multirow{2}{*}{\raisebox{2ex}{$0.797$}} & $\multirow{2}{*}{\raisebox{2ex}{\textbf{0.629}}}$ & $0.554$ & $ 0.482$ & $0.478$ & $0.476$ & $0.479$ & $\textbf{0.475}$ & $0.483$ & $0.558$\\
		&(0.107)& (0.255) & (0.192) & (0.135) & (0.091) & (0.069) & (0.062) & (0.060) & (0.060) & (0.063) & (0.060) & (0.060) & (0.056)\\
		\hline  
		\multirow{2}{*}{\ref{mod:4}} & \multirow{2}{*}{\raisebox{2ex}{$7.783$}} & \multirow{2}{*}{\raisebox{2ex}{$6.550$}} & \multirow{2}{*}{\raisebox{2ex}{$10.238$}} & \multirow{2}{*}{\raisebox{2ex}{$ 3.796$}} & $\multirow{2}{*}{\raisebox{2ex}{\bf 3.774}}$ & $3.608$ & $3.231$ & $3.185$ & $3.153$ & $3.189$ & $2.996$ & $3.186$ & $\bf 2.883$\\
		& (1.121) & (1.115) & (1.398) & (0.840) & (0.523) & (0.526) & (0.383) & (0.382) & (0.384) & (0.371) & (0.384) & (0.464) & (0.212)\\
		\hline  
		\multirow{2}{*}{\ref{mod:5}} & \multirow{2}{*}{\raisebox{2ex}{$0.508$}} & \multirow{2}{*}{\raisebox{2ex}{$0.518$}} & \multirow{2}{*}{\raisebox{2ex}{$0.699$}} & \multirow{2}{*}{\raisebox{2ex}{$0.575$}} & $\multirow{2}{*}{\raisebox{2ex}{\bf 0.436}}$ & $0.429$ & $0.389$ & $0.387$ & $0.386$ & $0.387$ & $\textbf{0.383}$ & $0.387$ & $0.486$\\
		&(0.051) & (0.073)& (0.084) & (0.081) & (0.051) &$(0.035)$ & $(0.031)$ & $(0.030)$ & $(0.030)$ & $(0.030)$ & $(0.030)$ & (0.028) & (0.077)\\
		\hline
		\multirow{2}{*}{\ref{mod:6}} & \multirow{2}{*}{\raisebox{2ex}{$\bf 1.015$}} & \multirow{2}{*}{\raisebox{2ex}{$1.020$}} & \multirow{2}{*}{\raisebox{2ex}{$1.405$}} & \multirow{2}{*}{\raisebox{2ex}{$1.774$}} & $\multirow{2}{*}{\raisebox{2ex}{1.290}}$ & $1.004$ & $0.934$ & $0.943$ & $0.941$ & $0.947$ & $\bf 0.914$ & $0.936$ & $0.957$\\
		&(0.054) & (0.053)& (0.098) & (0.145) & (0.083) &$( 0.085)$ & $(0.050)$ & $(0.062)$ & $(0.060)$ & $(0.053)$ & $(0.049)$ & (0.049) & (0.076)\\
		\hline
		\multirow{2}{*}{\ref{mod:7}} & \multirow{2}{*}{\raisebox{2ex}{$\bf 1.887$}} & \multirow{2}{*}{\raisebox{2ex}{$1.893$}} & \multirow{2}{*}{\raisebox{2ex}{$ 2.408$}} & \multirow{2}{*}{\raisebox{2ex}{$2.870$}} & $\multirow{2}{*}{\raisebox{2ex}{2.152}}$ & $1.939$ & $1.858$ & $1.854$ & $1.851$ & $1.867$ & $\bf 1.828$ & $1.852$ & $1.998$\\
		&(0.105) & (0.105)& (0.125) & (0.201) & (0.116) &$(0.109)$ & $(0.097)$ & $(0.097)$ & $(0.097)$ & $(0.098)$ & $(0.094)$ & (0.096) & (0.160)\\ 
		\hline  
		\multirow{2}{*}{\ref{mod:8}} & \multirow{2}{*}{\raisebox{2ex}{$1.475$}} & \multirow{2}{*}{\raisebox{2ex}{$\bf 1.461$}} & \multirow{2}{*}{\raisebox{2ex}{$1.578$}} & \multirow{2}{*}{\raisebox{2ex}{$1.919$}} & $\multirow{2}{*}{\raisebox{2ex}{1.464}}$ & $1.426$ & $1.416$ & $\textbf{1.416}$ & $\textbf{1.415}$ & $\textbf{1.419}$ & $\textbf{1.415}$ & $1.416$ & $1.456$\\
		&(0.079) & (0.078)& (0.089) & (0.121) & (0.074) &$(0.085)$ & $(0.080)$ & $(0.080)$ & $(0.080)$ & $(0.081)$ & $(0.079)$ & (0.080) & (0.099)\\
		\hline 
		\multirow{2}{*}{\ref{mod:9}} & \multirow{2}{*}{\raisebox{2ex}{$\bf 3.343$}} & \multirow{2}{*}{\raisebox{2ex}{$3.581$}} & \multirow{2}{*}{\raisebox{2ex}{$3.885$}} & \multirow{2}{*}{\raisebox{2ex}{$4.656$}} & $\multirow{2}{*}{\raisebox{2ex}{3.436}}$ & $3.332$ & $3.279$ & $3.279$ & $3.273$ & $3.293$ & $\textbf{3.240}$ & $3.277$ & $3.592$\\  
		&(0.187) & (0.499)  & (0.199)  & (0.292) & (0.186) & (0.172) & (0.164) & (0.170) & (0.169) & (0.175) & (0.167) & (0.168) & (0.176)\\
		\hline    
		\multirow{2}{*}{\ref{mod:10}} & \multirow{2}{*}{\raisebox{2ex}{$\bf 2.328$}} & \multirow{2}{*}{\raisebox{2ex}{$2.489$}} & \multirow{2}{*}{\raisebox{2ex}{$2.797$}} & \multirow{2}{*}{\raisebox{2ex}{$3.381$}} & $\multirow{2}{*}{\raisebox{2ex}{2.541}}$ & $2.308$ & $2.214$ &  $2.216$ & $2.212$ & $2.232$ & $\bf 2.171$ & $2.210$ & $2.342$\\
		&(0.135) & (0.158) & (0.154) & (0.243) & (0.141) & (0.163) & (0.143) & (0.154) &  (0.153)  &  (0.163) & (0.153) & (0.153) & (0.158)\\
		\hline                             
	\end{tabular}  
	\caption{Average MSEs in the correlated case.}
	\label{tab:corr}
	
	\hspace{0.5em}
	
	\begin{tabular}{ c | c c c c c | c c c c c c c c}  
		\hline                     
		\textcolor{cyan}{\bf Model} & \textcolor{cyan}{\bf Las} &\textcolor{cyan}{\bf Rid} &\textcolor{cyan}{\bf $k$NN} & \textcolor{cyan}{\bf Tr} & \textcolor{cyan}{\bf RF} & \textcolor{cyan}{\bf COBRA} & \textcolor{cyan}{\bf Epan} & \textcolor{cyan}{\bf Bi-wgt} & \textcolor{cyan}{\bf Tri-wgt}  & \textcolor{cyan}{\bf C-Gaus} & \textcolor{cyan}{\bf Gauss} & \textcolor{cyan}{\bf Exp$4$} & \textcolor{cyan}{\bf KCOBRA}\\ [0.5ex]   
		\hline  
		\multirow{2}{*}{\ref{mod:1}} & \multirow{2}{*}{\raisebox{2ex}{$2.294$}} & \multirow{2}{*}{\raisebox{2ex}{$1.947$}} & \multirow{2}{*}{\raisebox{2ex}{$1.941$}} & \multirow{2}{*}{\raisebox{2ex}{$\bf 0.320$}} & $\multirow{2}{*}{\raisebox{2ex}{0.542}}$ & $0.307$ & $0.304$ & $0.301$ & $0.288$ & $0.297$ & $\textbf{0.269}$ & $0.291$ & $0.449$\\
		&(0.544 & (0.507)& (0.487) & (0.145) & (0.231) &$(0.129)$ & $(0.105)$ & $(0.111)$ & $(0.103)$ & $(0.104)$ & $(0.092)$ & (0.098) & (2.50)\\
		\hline
		\multirow{2}{*}{\ref{mod:2}} & \multirow{2}{*}{\raisebox{2ex}{$14.273$}} & \multirow{2}{*}{\raisebox{2ex}{$8.442$}} & \multirow{2}{*}{\raisebox{2ex}{$8.572$}} & \multirow{2}{*}{\raisebox{2ex}{$ 6.796$}} & $\multirow{2}{*}{\raisebox{2ex}{\textbf{5.135}}}$ & $5.345$ & $4.582$ & $4.529$ & $4.491$ & $4.541$ & $\textbf{4.377}$ & $4.910$ & $4.946$\\
		&(2.593) & (1.912)& (1.751) & (1.548) & (1.372) &$(1.194)$ & $(0.941)$ & $(0.934)$ & $(0.922)$ & $(0.896)$ & $(0.905)$ & (1.181) & (1.271)\\
		\hline
		\multirow{2}{*}{\ref{mod:3}} & \multirow{2}{*}{\raisebox{2ex}{$7.996$}} & \multirow{2}{*}{\raisebox{2ex}{$6.266$}} & \multirow{2}{*}{\raisebox{2ex}{$8.704$}} & \multirow{2}{*}{\raisebox{2ex}{$4.110$}} & $\multirow{2}{*}{\raisebox{2ex}{\bf 3.722}}$ & $3.327$ & $2.598$ & $2.536$ & $2.444$ & $2.554$ & $2.168$ & $2.357$ & $\bf 1.853$\\
		&(3.393) & (3.296)& (3.523) & (2.894) & (2.956) & $(1.006)$ & $(0.912)$ & $(0.944)$ & $(0.840)$ & $(0.907)$ & $(0.680)$ & (0.756) & (0.443)\\
		\hline
		\multirow{2}{*}{\ref{mod:4}} & \multirow{2}{*}{\raisebox{2ex}{$61.474$}} & \multirow{2}{*}{\raisebox{2ex}{$42.351$}} & \multirow{2}{*}{\raisebox{2ex}{$46.934$}} & \multirow{2}{*}{\raisebox{2ex}{$\bf 8.855$}} & $\multirow{2}{*}{\raisebox{2ex}{13.381}}$ & $9.599$ & $10.511$ & $9.963$ & $9.682$ & $10.085$ & $9.056$ & $9.713$ & $\bf 8.957$\\
		&(13.986) & (11.622)& (12.543) & (3.480) & (5.549) &$(4.125)$ & $(2.961)$ & $(3.101)$ & $(2.860)$ & $(2.904)$ & $(2.407)$ & (2.695) & (0.954)\\
		\hline
		\multirow{2}{*}{\ref{mod:5}} & \multirow{2}{*}{\raisebox{2ex}{$6.805$}} & \multirow{2}{*}{\raisebox{2ex}{$7.479$}} & \multirow{2}{*}{\raisebox{2ex}{$10.342$}} & \multirow{2}{*}{\raisebox{2ex}{$\bf 4.000$}} & $\multirow{2}{*}{\raisebox{2ex}{4.880}}$ & $3.225$ & $2.640$ & $2.401$ & $2.235$ & $2.412$ & $\textbf{1.792}$ & $2.194$ & $2.873$\\
		&(3.685) & (5.336)& (5.425) & (3.144) & (3.787) &$(2.088)$ & $(1.455)$ & $(1.387)$ & $(1.250 )$ & $(1.355)$ & $(0.913)$ & (1.242) & (0.750)\\
		\hline   
		\multirow{2}{*}{\ref{mod:6}} & \multirow{2}{*}{\raisebox{2ex}{$24.078$}} & \multirow{2}{*}{\raisebox{2ex}{$23.883$}} & \multirow{2}{*}{\raisebox{2ex}{$22.216$}} & \multirow{2}{*}{\raisebox{2ex}{$24.612$}} & $\multirow{2}{*}{\raisebox{2ex}{\bf 20.202}}$ & $19.573$ & $18.475$ & $18.901$ & $16.718$ & $17.186$ & $\textbf{14.982}$ & $16.597$ & $18.541$\\
		&(5.547) & (5.527) & (5.255) & (5.351) & (5.291) & (5.919) & (4.886) & (5.703) & (5.569) & (6.232) & (5.556) & (5.479) & (6.863)\\
		\hline   
		\multirow{2}{*}{\ref{mod:7}} & \multirow{2}{*}{\raisebox{2ex}{$2.358$}} & \multirow{2}{*}{\raisebox{2ex}{$2.357$}} & \multirow{2}{*}{\raisebox{2ex}{$2.602$}} & \multirow{2}{*}{\raisebox{2ex}{$2.890$}} & $\multirow{2}{*}{\raisebox{2ex}{\bf 2.260}}$ & $2.312$ & $2.221$ & $2.223$ & $2.220$ & $2.236$ & $\textbf{2.216}$ & $2.223$ & $2.294$\\
		&(0.122) & (0.122) & (0.125) & (0.165) & (0.112) & (0.124) & (0.106) & (0.112) & (0.111) & (0.121) & (0.111) & (0.110) & (0.168)\\
		\hline    
		\multirow{2}{*}{\ref{mod:8}} & \multirow{2}{*}{\raisebox{2ex}{$4.013$}} & \multirow{2}{*}{\raisebox{2ex}{$\bf 3.929$}} & \multirow{2}{*}{\raisebox{2ex}{$ 4.276$}} & \multirow{2}{*}{\raisebox{2ex}{$5.151$}} & $\multirow{2}{*}{\raisebox{2ex}{3.986}}$ & $4.046$ & $3.948$ & $3.953$ & $3.949$ & $3.973$ & $\textbf{3.937}$ & $\bf 3.948$ & $4.213$\\
		&(0.258) & (0.253) & (0.259) & (0.394) & (0.241) & (0.253) & (0.225) & (0.244) & (0.246) & (0.249) & (0.244) & (0.245) & (0.341)\\
		\hline
		\multirow{2}{*}{\ref{mod:9}} & \multirow{2}{*}{\raisebox{2ex}{$\bf 6.072$}} & \multirow{2}{*}{\raisebox{2ex}{$9.764$}} & \multirow{2}{*}{\raisebox{2ex}{$8.308$}} & \multirow{2}{*}{\raisebox{2ex}{$10.647$}} & $\multirow{2}{*}{\raisebox{2ex}{7.862}}$ & $6.450$ & $6.017$ & $5.954$ & $5.906$ & $5.923$ & $\bf 5.732$ & $5.841$ & $6.407$\\
		&(0.672) & (0.610) & (0.593) & (0.645) & (0.560) & (0.629) & (0.527) & (0.572) & (0.566) & (0.516) & (0.498) & (0.515) & (0.950)\\
		\hline  
		\multirow{2}{*}{\ref{mod:10}} & \multirow{2}{*}{\raisebox{2ex}{$\bf 15.402$}} & \multirow{2}{*}{\raisebox{2ex}{$17.611$}} & \multirow{2}{*}{\raisebox{2ex}{$19.287$}} & \multirow{2}{*}{\raisebox{2ex}{$20.819$}} & $\multirow{2}{*}{\raisebox{2ex}{16.039}}$ & $15.754$ & $16.629$ & $14.970$ & $15.017$ & $14.806$ & $\textbf{14.346}$ & $14.568$ & $16.666$\\
		&(1.644) & (2.225) & (1.885) & (1.844) & (1.748) & (1.825) & (2.348) & (1.913) & (1.967) & (1.705) & (1.506) & (1.517) & (1.534)\\
		\hline        
	\end{tabular}
\end{sidewaystable}

\begin{figure}[h!]
	\includegraphics[width=7.4cm]{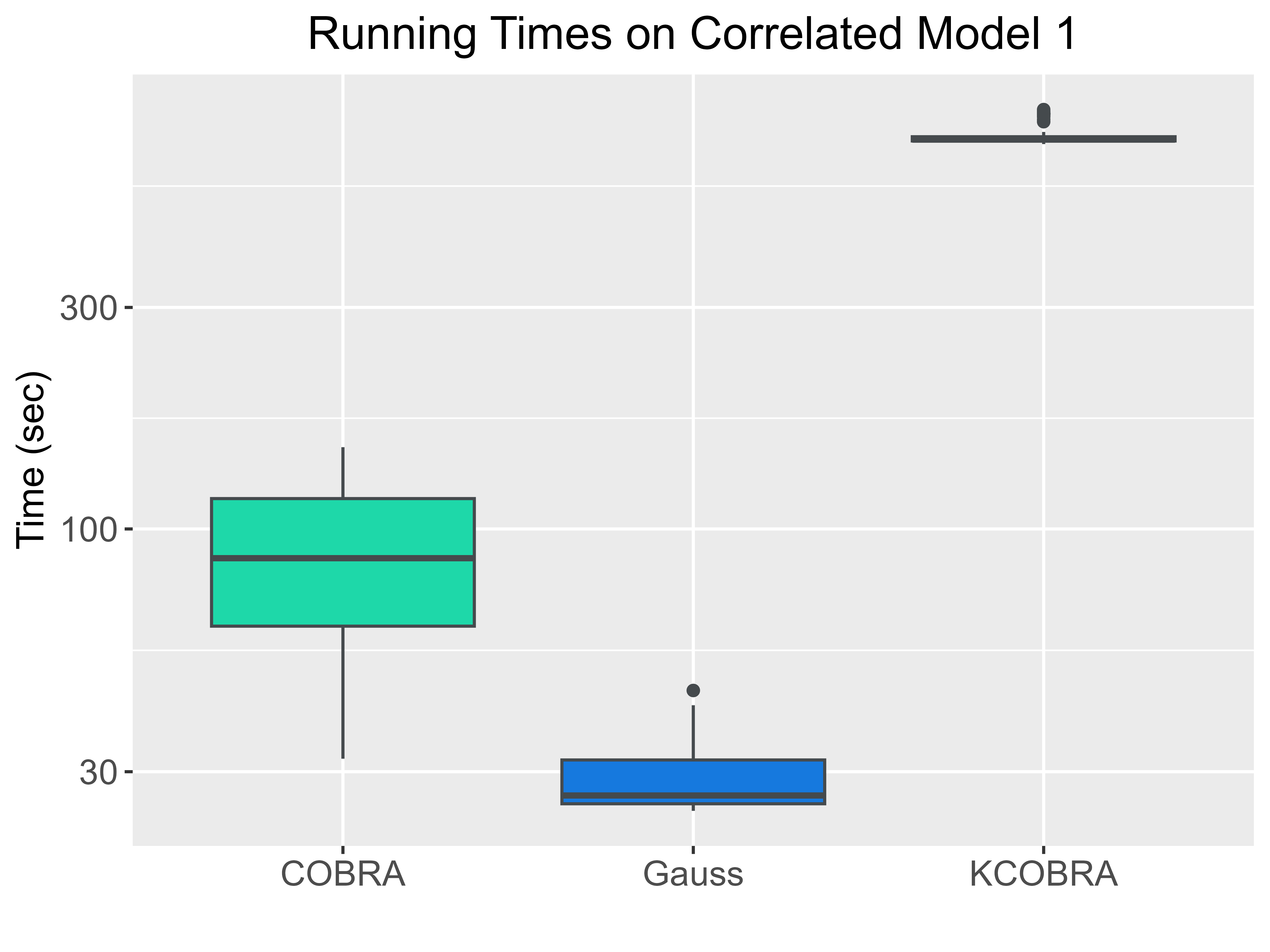}\noindent
	\includegraphics[width=7.4cm]{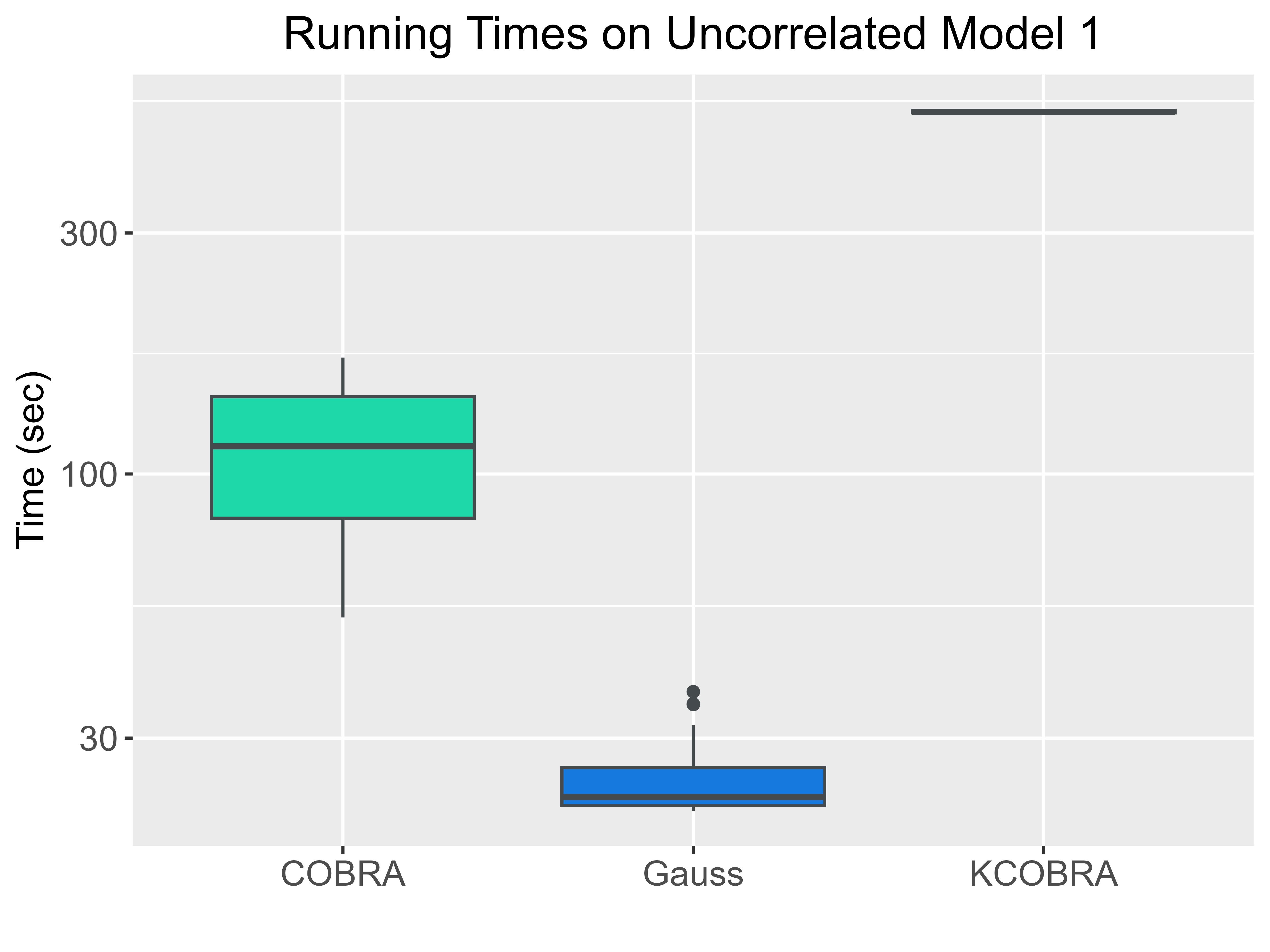}
	\includegraphics[width=7.4cm]{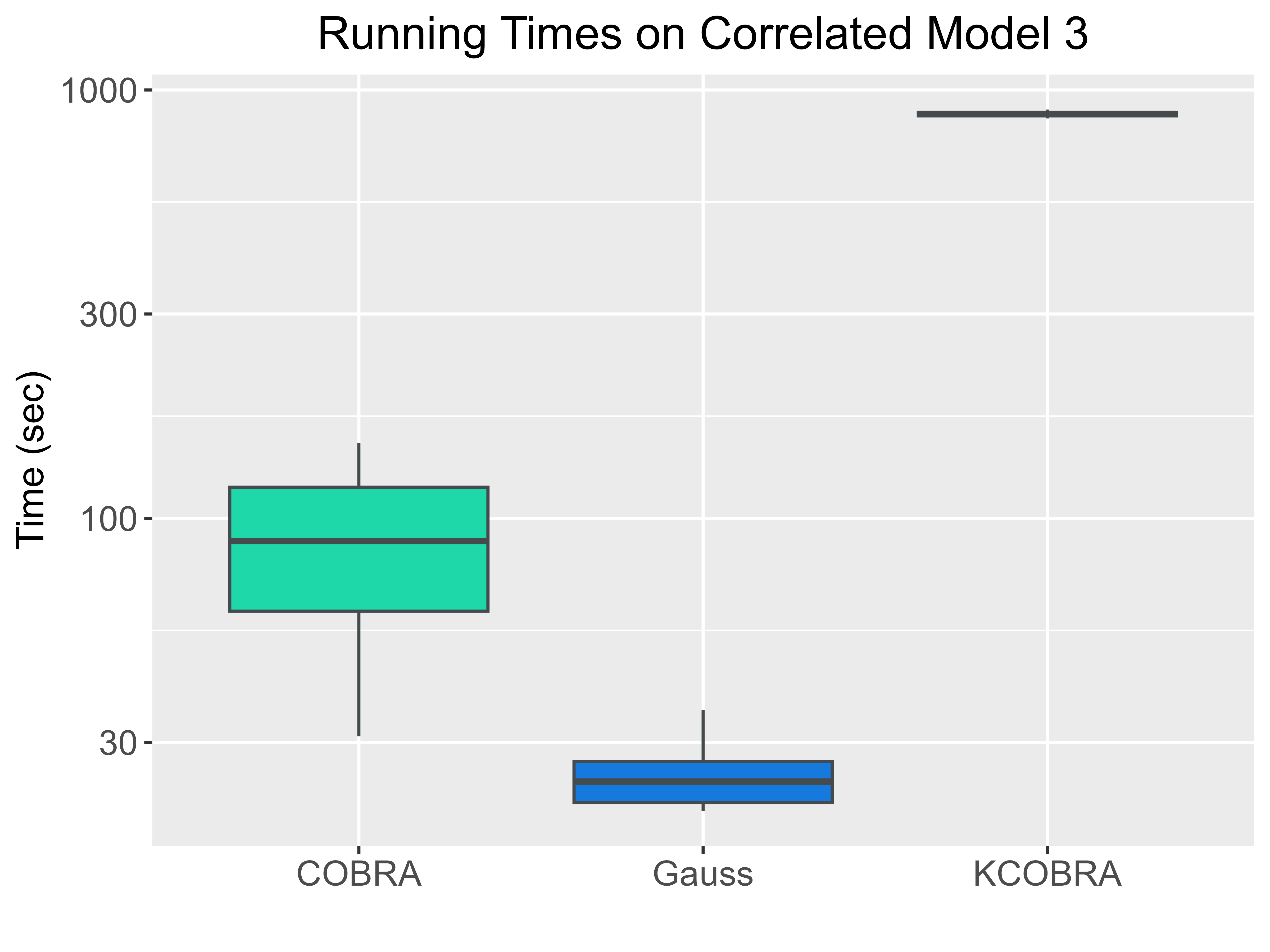}\noindent
	\includegraphics[width=7.4cm]{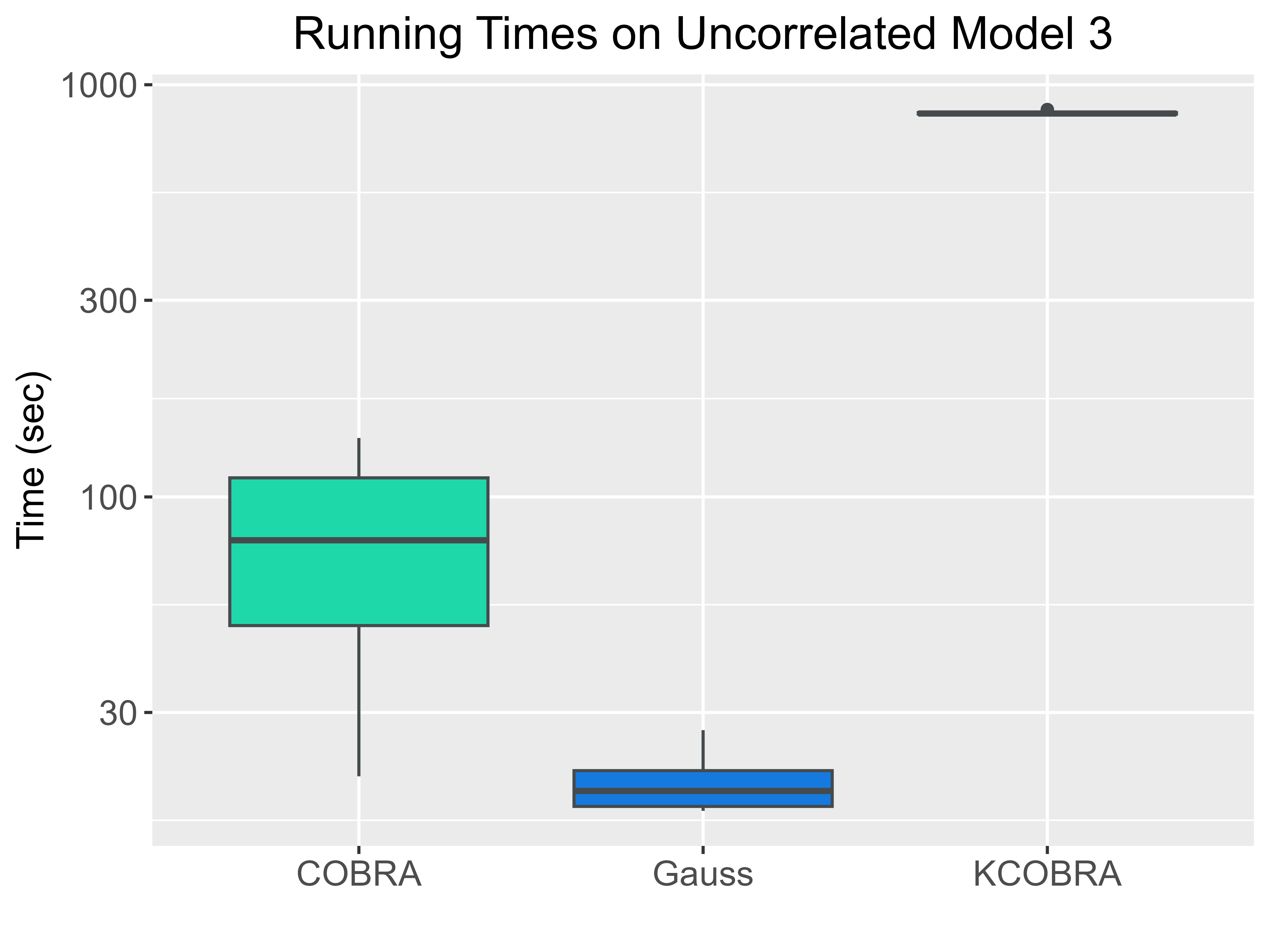}
	\caption{Boxplots of computational times of the three aggregation strategies implemented on model~\ref{mod:1} and \ref{mod:3}, with $500$ of bandwidth parameters. Note that ``Gauss'' corresponds to the proposed method with Gaussian kernel, and the ``Time'' axis is in logarithmic scale.}
	\label{fig:2}
\end{figure}

In each table, we are interested in comparing the smallest average RMSEs in the first block to all the columns in the second block. First of all, we can see that all columns of the second block often outperform the best estimator of the first block, which illustrates the theoretical result of the combining estimation methods. Secondly, the proposed methods (second to seventh column of the second block) always outperform the classical COBRA (first column) and KernelCOBRA (last column) for almost all kernels. Lastly, the combining estimation method with Gaussian kernel is the best one in both tables. In addition, Figure~\ref{fig:2} below contains boxplots of RMSEs obtained from 100 independent runs of Model \ref{mod:1} and \ref{mod:10} (correlated and uncorrelated cases), computed on a computational machine with the following characteristics:

\begin{itemize}
	\item Processor: 2x AMD Opteron 6174, 12C, 2.2GHz, 12x512K L2/12M L3 Cache, 80W ACP, DDR3-1333MHz.
	\item Memory: 64GB Memory for 2 CPUs, DDR3, 1333MHz.
\end{itemize}

\noindent These boxplots clearly show that the proposed method is around 3 to 10 times faster than the classical method by \cite{cobra}, and is up to hundred times faster than KernelCOBRA by \cite{pycobra} with $500$ values of bandwidth parameters.

\subsection{Real public datasets}
In this part, we consider three public datasets which are available and easily accessible on the internet. The first dataset ({\bf Abalone}, available at \cite{abaloneData}) contains $4177$ rows and $9$ columns of measurements of abalones observed in Tasmania, Australia. We are interested in predicting the age of each abalone through the number of rings using its physical characteristics such as gender, size, weight,  etc. The second dataset ({\bf House}, available at \cite{HouseKC2016}) comprises house sale prices for King County including Seattle. It contains homes sold between May 2014 and May 2015. The dataset consists of $21613$ rows of houses and $21$ columns of characteristics of each house including ID, Year of sale, Size, Location, etc. In this case, we want to predict the price of each house using all of its quantitative characteristics. 

Finally, the last dataset ({\bf Wine}, see \cite{redWineData,wineArticle}), which was also considered in \cite{cobra}, containing $1599$ rows of different types of wines, and $12$ columns corresponding to different substances of red wines including the amount of different types of acids, sugar, chlorides, PH, etc. The variable of interest is {\it quality} which scales from $3$ to $8$ where $8$ represents the best quality. We aim at predicting the quality of each wine, which is treated as a continuous variable, using all of its substances.

The five primary regressors are Ridge, LASSO, $k$NN, Tree and Random Forest regression. In this case, the parameter $ntree=500$ for random forest, and $k$NN is implemented using $k=20,12$ and $5$ for Abalone, House and Wine dataset respectively. The five regressors are combined using the classical method by \cite{cobra}, the proposed method using Gaussian kernel, and the KernelCOBRA by \cite{pycobra}. In this case, $300$ values of parameter $h$ are considered for the classical COBRA and KernelCOBRA. 

The average RMSEs obtained from $100$ independent runs, evaluated on $20\%$-testing data of the three public datasets, are provided in Table~\ref{tab:realdata} below (the first three rows). We observe that random forest is the best estimator among all the basic estimators in the first block, and the proposed method (\textcolor{cyan}{\bf Gauss}) either outperforms other columns ({\bf Wine} and {\bf Abalone}) or biases towards the best basic estimator ({\bf House}). Moreover, the performances of the proposed method always exceed the ones of the classical COBRA and the KernelCOBRA.

\subsection{Real private datasets}
In this section, we provide the performances of the aggregation methods on other two (real) private datasets. The first dataset contains six columns corresponding to the six variables including \textit{Air temperature, Input Pressure, Output Pressure, Flow, Water Temperature} and \textit{Power Consumption} along with $2026$ rows of hourly observations of these measurements of an air compressor machine provided by \cite{CHM}. The goal is to predict the power consumption of this machine using the five remaining explanatory variables. The second dataset is provided by the wind energy company Ma$\ddot{\text{\i}}$a Eolis. It contains 8721 observations of seven variables representing 10-minute measurements of \textit{Electrical power, Wind speed, Wind direction, Temperature, Variance of wind speed} and \textit{Variance of wind direction} measured from a wind turbine of the company (see, \cite{wind}). In this case, we aim at predicting the electrical power produced by the turbine using the remaining six measurements as explanatory variables. We use the same set of parameters as in the previous subsection except for $k$NN where in this case $k=10$ and $k=7$ are used for air compressor and wind turbine dataset respectively. 

\begin{table}[!ph]
	\renewcommand{\arraystretch}{1.3}
	\centering
	\tiny
	\caption{Average RMSEs of real datasets.}
	\label{tab:realdata}
	\begin{tabular}{l@{\hspace{1em}} c @{\hspace{1em}} c @{\hspace{1em}} c @{\hspace{1em}} c @{\hspace{1em}} c @{\hspace{1em}} | c @{\hspace{1em}} c @{\hspace{1em}} c}  
		\hline                     
		\textcolor{cyan}{\bf Model} & \textcolor{cyan}{\bf Las} &\textcolor{cyan}{\bf Rid} &\textcolor{cyan}{\bf $k$NN} & \textcolor{cyan}{\bf Tr} & \textcolor{cyan}{\bf RF} & \textcolor{cyan}{\bf COBRA} & \textcolor{cyan}{\bf Gauss} & \textcolor{cyan}{\bf KCOBRA}\\ [0.5ex] 
		\hline  
		\multirow{2}{*}{\bf Abalone} & \multirow{2}{*}{\raisebox{2ex}{2.20}} & \multirow{2}{*}{\raisebox{2ex}{2.22}} & \multirow{2}{*}{\raisebox{2ex}{2.18}} & \multirow{2}{*}{\raisebox{2ex}{2.40}} & \multirow{2}{*}{\raisebox{2ex}{\bf 2.15}} & \multirow{2}{*}{\raisebox{2ex}{2.17}} & \textbf{2.13} & 2.67\\
		& (0.07) & (0.08)& (0.06) & (0.07) & (0.06) &(0.08) & (0.06) & (0.12)\\
		\hline
		\multirow{2}{*}{\bf House} & \multirow{2}{*}{\raisebox{2ex}{241083.96}} & \multirow{2}{*}{\raisebox{2ex}{241072.97}} & \multirow{2}{*}{\raisebox{2ex}{245153.61}} & \multirow{2}{*}{\raisebox{2ex}{254099.65}} & \multirow{2}{*}{\raisebox{2ex}{\bf 205943.77}} & \multirow{2}{*}{\raisebox{2ex}{223596.32}} & \textbf{209955.28} & {650943.60} \\
		& (8883.11) & (8906.33)& (23548.37) & (9350.89) & (7496.77) & (13299.93) & (7815.62) & (29565.23)\\
		\hline
		\multirow{2}{*}{\bf Wine} & \multirow{2}{*}{\raisebox{2ex}{0.66}} & \multirow{2}{*}{\raisebox{2ex}{0.69}} & \multirow{2}{*}{\raisebox{2ex}{0.77}} & \multirow{2}{*}{\raisebox{2ex}{0.71}} & \multirow{2}{*}{\raisebox{2ex}{\bf 0.62}} & 0.65 & $\multirow{2}{*}{\raisebox{2ex}{\textbf{0.62}}}$ & 0.67\\
		& (0.03) & (0.05)& (0.03) & (0.03) & (0.03) & (0.03) & (0.02) & (0.04)\\
		\hline
		\multirow{2}{*}{\bf Air} & \multirow{2}{*}{\raisebox{2ex}{\bf163.10}} & \multirow{2}{*}{\raisebox{2ex}{164.23}} & \multirow{2}{*}{\raisebox{2ex}{241.66}} & \multirow{2}{*}{\raisebox{2ex}{351.32}} & \multirow{2}{*}{\raisebox{2ex}{174.84}} & $\multirow{2}{*}{\raisebox{2ex}{172.86}}$ & \textbf{163.25} & 1468.30\\
		& (3.69) & (3.75)& (5.87) & (31.88) & (6.55) &$(7.64)$ & $(3.33)$ & (78.47)\\
		\hline
		\multirow{2}{*}{\bf Turbine} & \multirow{2}{*}{\raisebox{2ex}{70.05}} & \multirow{2}{*}{\raisebox{2ex}{68.99}} & \multirow{2}{*}{\raisebox{2ex}{44.52}} & \multirow{2}{*}{\raisebox{2ex}{81.71}} & \multirow{2}{*}{\raisebox{2ex}{$\bf 38.89$}} & $\multirow{2}{*}{\raisebox{2ex}{38.93}}$ & \textbf{37.14} & 515.41\\
		& (4.99) & (3.41)& (1.67) & (4.98) & (1.51) &$(1.56)$ & $(1.56)$ & (58.14)\\
		\hline
	\end{tabular}
\end{table}

The results obtained from $100$ independent runs of the methods are presented in the last two rows (\textbf{Air} and \textbf{Turbine}) of Table~\ref{tab:realdata} above. We observe on one hand that the proposed method (\textcolor{cyan}{\bf Gauss}) outperforms the best basic estimators (\textcolor{cyan}{\bf RF}) and the two other competitors (\textcolor{cyan}{\bf COBRA} and \textcolor{cyan}{\bf KCOBRA}) in the case of {\bf Turbine} dataset. On the other hand, the performance of our method approaches the best basic estimator (\textcolor{cyan}{\bf Las}) and outperforms the other aggregation methods in the case of \textbf{Air} dataset.

Moreover, boxplots of 100 runs measured on {\bf Wine} and {\bf Turbine} datasets (computed using the same computational machine as described in section~\ref{subsec:simulated}) are also provided in Figure~\ref{fig:3} below.

\begin{figure}[!h]
	\includegraphics[width=7cm]{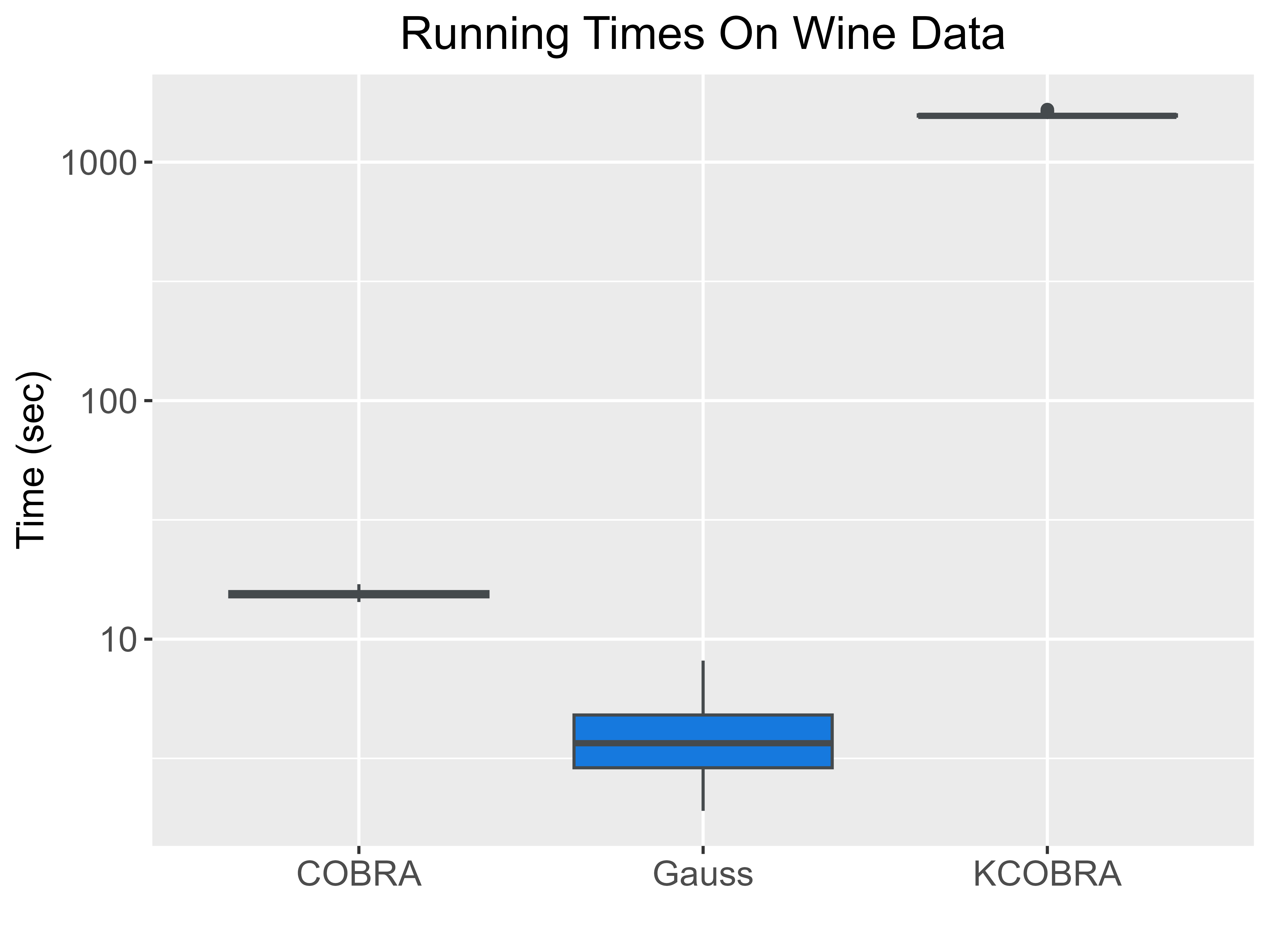}
	\includegraphics[width=7cm]{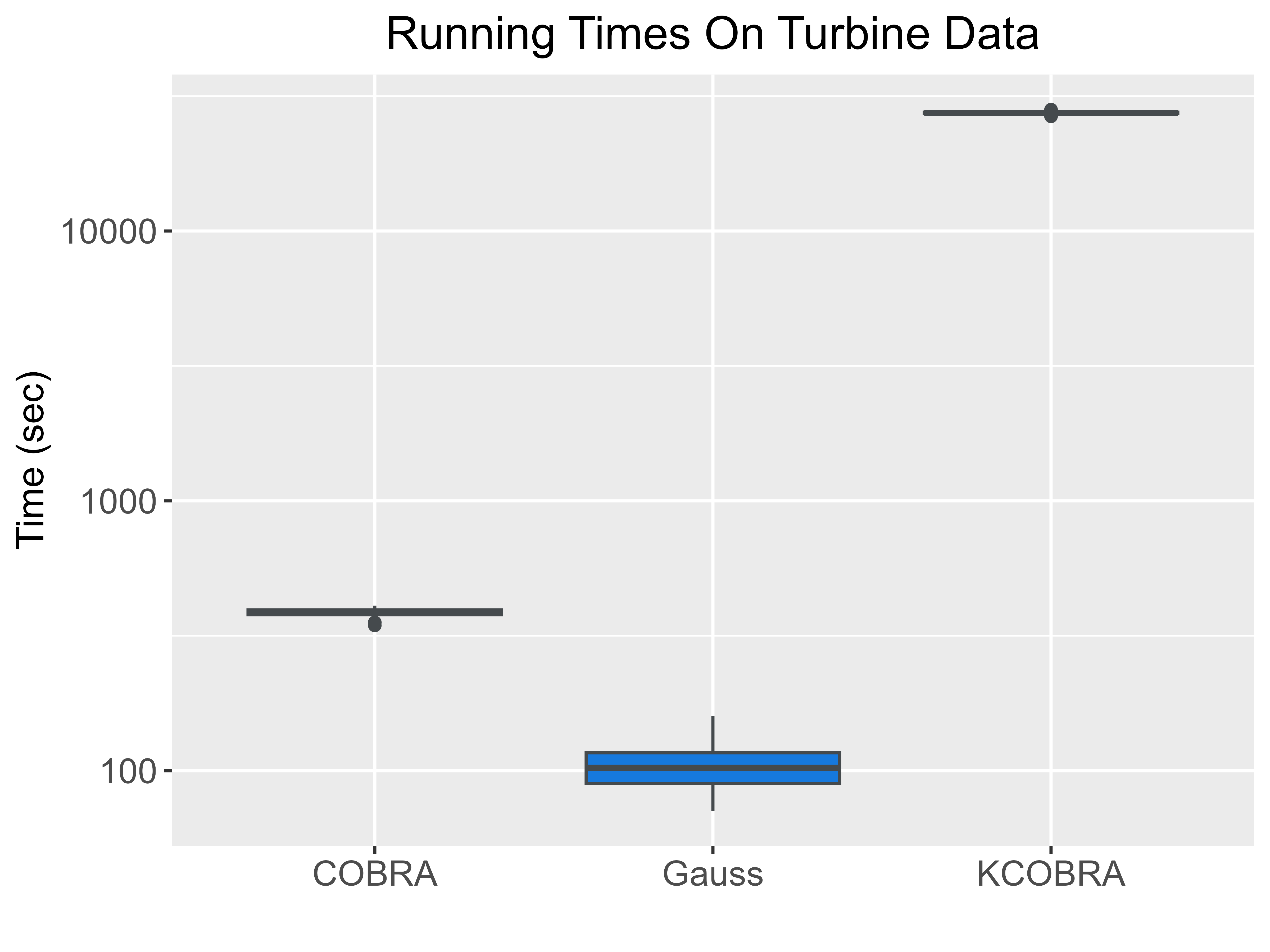}
	\caption{Boxplots of computational times of the proposed method and the two competitors implemented on {\bf Wine} and {\bf Turbine} datasets.}
	\label{fig:3}
\end{figure}

\subsection{Application on a data of Magnetosphere- Ionosphere System provided by CEA}
\label{sec:applicationCEA}
This section presents an application of the proposed method on a data provided by researchers of Commissariat à l'Énergie Atomique (CEA). In a collaboration with researchers of CEA on a research topic in Magnetosphere-Ionosphere System (see \cite{MLPitchAngle2022}), we are interested in constructing a global machine learning model of event-driven for estimating a physical quantity called {\it Pitch Angle Diffusion Coefficient} ($D_{\alpha\alpha}$) using three input data: electron at L-shell $L$, energy $E$, and equatorial pitch angle $\alpha$. Pitch angle diffusion coefficient is one of the major mechanisms that drives the structure of the Van Allen radiation belts and causes the well-known two belt structure. Whistler mode waves which are known to play a crucial role in thermodynamics, electron acceleration, and electron precipitation in the atmosphere are also caused by the physical process of pitch angle diffusion. This quantity can be computed from statistical models derived from years of satellite observations of the hiss waves properties of different missions, or using a method called event-driven approach (\cite{Thor2013}). We use in this study a database of event-driven diffusion coefficients that was generated for the studies of \cite{Ripo2019}. A very large fully observed dataset containing around two hundred million observations is available. However, one wants to construct predictive models using reasonably small training data, therefore, a 3-dimensional grid made up of $4$ values of $L\in\{2,3,4,5\}$, $60$ values of $E$ and $256$ values of $\alpha$ is considered. This filtering process creates a training dataset of size $61\ 440$, simply called $\mathcal{D}_0$. Then, two training datasets are extracted: high-resolution ($\mathcal{D}_{\text{HR}}$) and low-resolution datasets ($\mathcal{D}_{\text{LR}}$). High-resolution dataset is composed of $84$ pitch angles ($\alpha$) and $60$ energies bins ($E$), thus contains $20\ 160$ data points. The low-resolution dataset is composed of only $14$ pitch angles and $13$ energies bins, thus contains only $728$ data points. The table~\ref{tab:dataCEA1} below provides the structure of the described training datasets.
\begin{table}[!h]
	\centering 
	\small
	\begin{tabular}{| l | c | c | c | l | }  
		\hline                      
		\textbf{Data} & \textbf{$L$} & \textbf{$E$} & \textbf{$\alpha$} & \textbf{Size}\\
		\hline
		$\mathcal{D}_{0}$ & $4$ & $60$ & $256$ & $61\ 440$\\ 
		\hline
		$\mathcal{D}_{\text{HR}}$ & $4$ & $60$ & $84$ & $20\ 160$\\ 
		\hline
		$\mathcal{D}_{\text{LR}}$ & $4$ & $13$ & $14$ & $728$\\ 
		\hline
	\end{tabular}
	\caption{The high and low resolution training datasets.}%
	\label{tab:dataCEA1}
\end{table}

\begin{figure}[!h]
	\centering
	\includegraphics[scale=0.7]{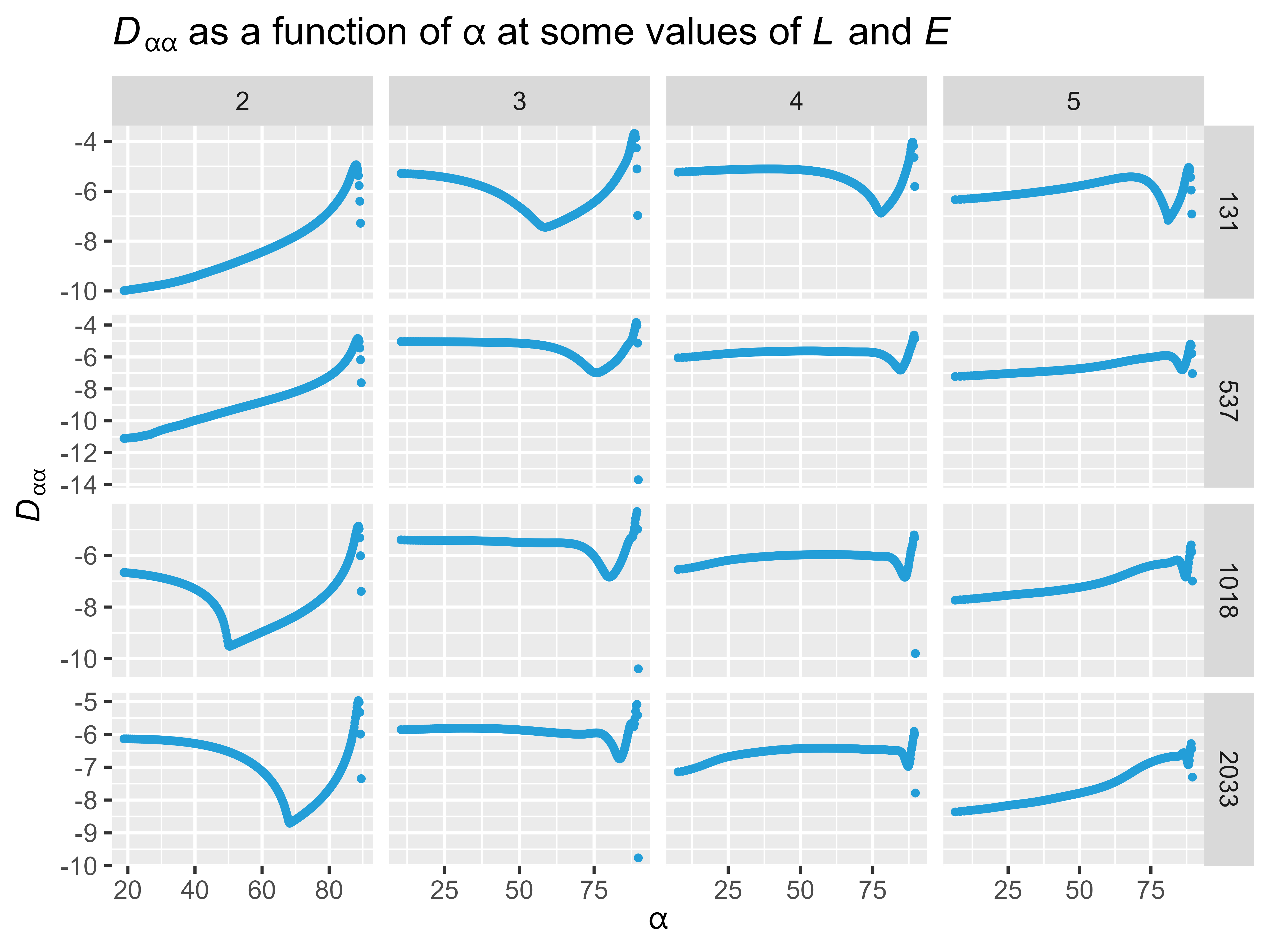}
	\caption{The relation between $D_{\alpha\alpha}$ and $\alpha$ at some cuts of $L$ and $E$ values.}
	\label{fig:D_alphaVSalpha}
\end{figure}

It should be pointed out that the training datasets are noiseless (see Figure~\ref{fig:D_alphaVSalpha}), and the relationship of $D_{\alpha\alpha}$ and $\alpha$ at some fixed couples $(L,E)$ are illustrated in Figure~\ref{fig:D_alphaVSalpha} below.

In this part, we considered several regression models, including $k$-nearest neighbors (kNN), kernel regression (KerReg), regression tree (Tree), bagging (Bag), random forest (RF), radial basis (Radial), splines (Spline), and deep neural networks (DNN). These models were trained separately on the high-resolution ($\mathcal{D}_{\text{HR}}$) and low-resolution ($\mathcal{D}_{\text{LR}}$) training datasets.

To evaluate the prediction capability of these models, we extracted three different testing datasets from the fully observed data, which contains two hundred million observations. By using these testing datasets, we were able to compare the performance of the different regression models and identify the most effective one for the task at hand. Table~\ref{tab:testCEA1} below describes the three testing datasets. 
\begin{table}[h!]
	\centering 
	\small
	\begin{tabular}{| l | l | }  
		\hline 
		\begin{tabular}{l}                     
			\textbf{Data}
		\end{tabular} & 
		\begin{tabular}{l}
			\textbf{Description}
		\end{tabular}\\
		\hline
		\begin{tabular}{l}
			${\cal D}_{\text{testHR}}$
		\end{tabular} & 
		\begin{tabular}{l}
			For testing the models built on $D_{\text{HR}}$.
		\end{tabular}\\ 
		\hline
		\begin{tabular}{l}
			${\cal D}_{\text{testLR}}$
		\end{tabular} & 
		\begin{tabular}{l}
			For testing the models built on $D_{\text{LR}}$.
		\end{tabular}\\
		\hline
		\begin{tabular}{l}
			${\cal D}_{\text{testL}}$
		\end{tabular} & 
		\begin{tabular}{l}Contains more values of $L$ other than $\{2,3,4,5\}$.\\ For testing the models built on both training data. \end{tabular}\\ 
		\hline
	\end{tabular}
	\caption{The three testing datasets.}%
	\label{tab:testCEA1}
\end{table}

In both cases, the regression estimators were constructed using the entire training data ($\mathcal{D}_{\text{HR}}$ or $\mathcal{D}_{\text{LR}}$), which left no training data for aggregation. To avoid violating the independence assumption between the data used to train the individual estimators and the data used for aggregation, we randomly divided each testing dataset into two parts. The first part is used to optimize the bandwidth parameter $h$ for the aggregation, while the remaining part is used as the actual testing dataset. The numerical results obtained from 50 independent runs of this procedure are presented in Figure~\ref{fig:RMSE_CEA} below.

\begin{figure}[h!]
	\centering
	\includegraphics[width=0.475\textwidth]{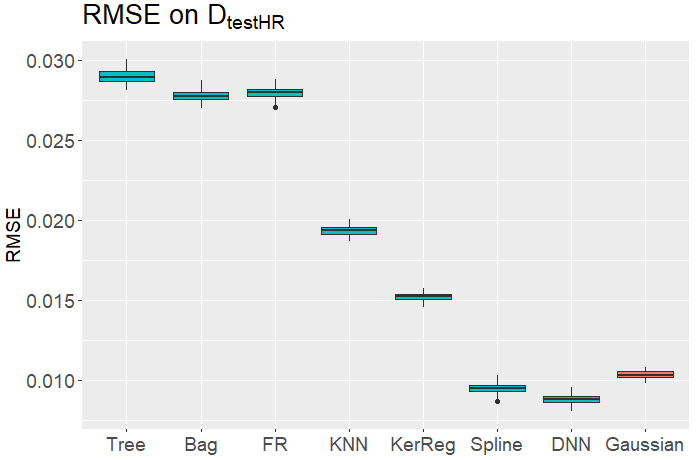}
	\includegraphics[width=0.475\textwidth]{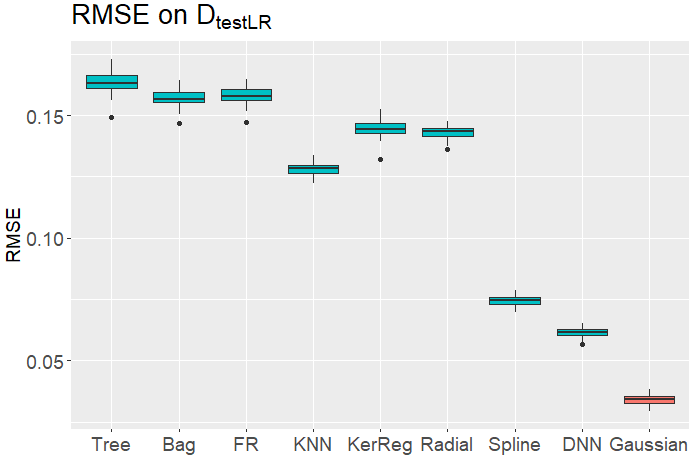}\\[0.3cm]
	\includegraphics[width=0.475\textwidth]{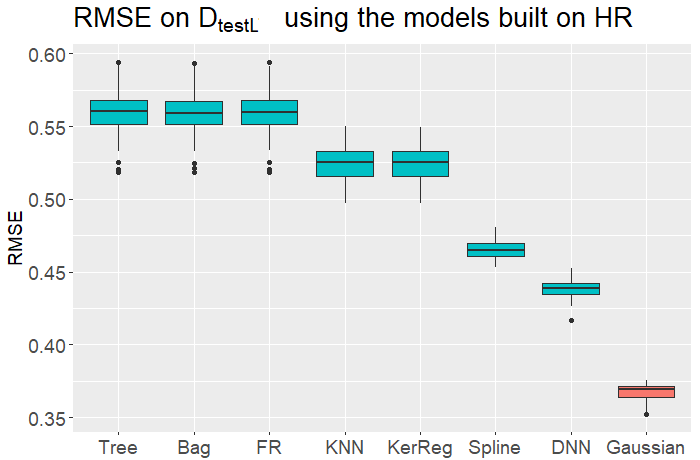}
	\includegraphics[width=0.475\textwidth]{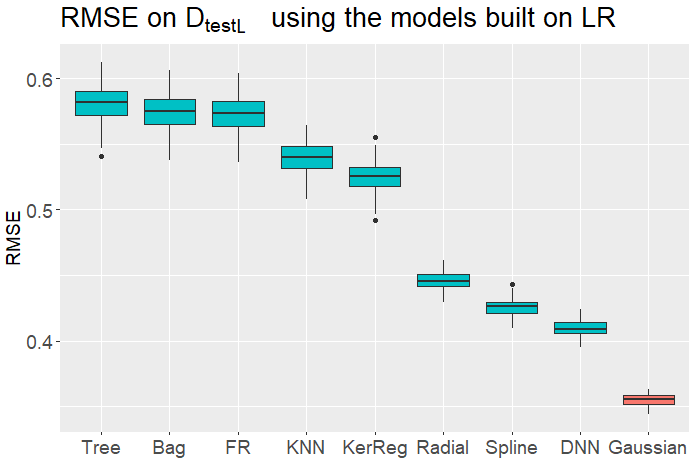}
	\caption{Boxplots of RMSEs over $50$ runs of the algorithm. Note that Radial is built only on the training data $\mathcal{D}_{\text{LR}}$, therefore it is not presented in the two boxplots on the left-hand side (where the model are built using $\mathcal{D}_{\text{HR}}$). The last boxplot is the performance of the proposed aggregation method using Gaussian kernel.}
	\label{fig:RMSE_CEA}
\end{figure}

The kernel-based consensual aggregation method is implemented using Gaussian kernel and is denoted by Gaussian. We observe that the tree-based models behave similarly and are the weak ones, and DNN is the best individual estimator as it provides the lowest average testing RMSE. On the other hand, the aggregation outperforms other basic estimators in the last three cases, and biases towards the best basic estimator on $D_{\text{testHR}}$.
\begin{remark}
	As the training data in our study are extracted selectively from the full observed data, the distributions of the training and testing data are not the same. For instance, $L$ only takes values in $\{2,3,4,5\}$ in the training data, while the testing data may have more decimal values. To overcome this limitation, we split the testing data into two parts, allowing us to fine-tune the smoothing parameter $h$ and adjust the weights for predicting new observations coming from a different distribution.
	
	This approach is practically useful because the basic models can be built on one source of an underlying distribution ${\cal L}_0$ and then used to predict observations from another source of distribution ${\cal L}_1$, which may be different from ${\cal L}_0$. In such cases, access to a part of the new source is required to adjust the weights in the aggregation, akin to a domain adaptation-like property. This adaptability of the aggregation method is a remarkable advantage and can lead to improved performance in diverse settings.
\end{remark}

\section{Conclusion}
\label{sec:Conclude}
In conclusion, this study extends the context of a naive kernel-based consensual regression aggregation method to a more general regular kernel-based framework, and it demonstrates the consistency inheritance property of the method with the same convergence rate. Additionally, we propose an optimization algorithm based on gradient descent to efficiently estimate the key parameter of the method with the computational speed up to several hundred times faster than the classical grid search. Our numerical simulations show that the performance of the method is significantly improved with smoother kernel functions. Furthermore, we demonstrate, in a real-world project with physics data, that the method exhibits a domain adaptation-like property, which opens up interesting directions for further study.

In practice, the performance of the consensual aggregation method depends on both the individual regression estimators and the final combination, which involves kernel functions. Therefore, calibration of hyperparameters in both steps is critical, and automated machine learning models may be useful for improving the performance of the global model.

\section{Reproducibility of the experiments}
\label{sec:supp}


For readers interested in reproducing our experiments, we have made some public datasets used in this article and the official source codes written in \texttt{python} and \texttt{R} of the algorithm available on our Github repository: \url{https://github.com/hassothea/AggregationMethods}.

\section{Proofs}
\label{sec:proof}

The following lemma, which is a variant of lemma 4.1 in \cite{bookDistributionFree} related to the property of binomial random variables, is needed. 

\begin{lemma}
	\label{lem:1}
	Let $B(n,p)$ be the binomial random variable with parameters $n$ and $p$. Then
	\begin{enumerate}
		\item\label{itm:1} For any $c>0$,
		\begin{align*}
			\mathbb{E}\Big[\frac{1}{c+B(n,p)}\Big]\leq \frac{2}{p(n+1)}.
		\end{align*}
		\item\label{itm:2} \begin{align*}
			\mathbb{E}\Big[\frac{1}{B(n,p)}\mathds{1}_{B(n,p)>0}\Big]\leq \frac{2}{p(n+1)}.
		\end{align*}
	\end{enumerate}
\end{lemma}

\begin{prooflemma}
	\begin{enumerate}
		\item For any $c>0$, one has
		\begin{align*}
			\mathbb{E}\Big[\frac{1}{c+B(n,p)}\Big]&=\sum_{k=0}^n\frac{1}{c+k}\times\frac{n!}{(n-k)!k!}p^k(1-p)^{n-k}\\
			&= \sum_{k=0}^n\frac{1}{k+1}\times\frac{k+1}{k+c}\times\frac{n!}{(n-k)!k!}p^k(1-p)^{n-k}\\
			&\leq\frac{2}{p(n+1)} \sum_{k=0}^{n}\frac{(n+1)!p^{k+1}(1-p)^{n+1-(k+1)}}{[n+1-(k+1)]!(k+1)!}\\
			&\leq\frac{2}{p(n+1)} \sum_{k=0}^{n+1}\frac{(n+1)!p^k(1-p)^{n+1-k}}{[n+1-k]!k!}\\
			&=\frac{2}{p(n+1)} (p+1-p)^{n+1}\\
			&=\frac{2}{p(n+1)}
		\end{align*}
		\item 
		\begin{align*}
			\mathbb{E}\Big[\frac{1}{B(n,p)}\mathds{1}_{B(n,p)>0}\Big]&\leq \mathbb{E}\Big[\frac{2}{1+B(n,p)}\Big]\\
			&= \sum_{k=0}^n\frac{2}{k+1}\times\frac{n!}{(n-k)!k!}p^k(1-p)^{n-k}\\
			&=\frac{2}{p(n+1)} \sum_{k=0}^{n}\frac{(n+1)!p^{k+1}(1-p)^{n+1-(k+1)}}{[n+1-(k+1)]!(k+1)!}\\
			&\leq\frac{2}{p(n+1)} \sum_{k=0}^{n+1}\frac{(n+1)!p^k(1-p)^{n+1-k}}{[n+1-k]!k!}\\
			&=\frac{2}{p(n+1)} (p+1-p)^{n+1}\\
			&=\frac{2}{p(n+1)}
		\end{align*}
	\end{enumerate}
	\QED
\end{prooflemma}


\begin{proofprop}
	For any square integrable function with respect to $\textbf{r}_k(X)$, one has
	
	\begin{align*}
		\mathbb{E}\Big[|g_n(\textbf{r}_k(X))-g^*(X)|^2\Big]&= \mathbb{E}\Big[|g_n(\textbf{r}_k(X))-g^*(\textbf{r}_k(X))+g^*(\textbf{r}_k(X))-g^*(X)|^2\Big]\\
		&=\mathbb{E}\Big[|g_n(\textbf{r}_k(X))-g^*(\textbf{r}_k(X))|^2\Big]\\
		&\quad+2\mathbb{E}\Big[(g_n(\textbf{r}_k(X))-g^*(\textbf{r}_k(X)))(g^*(\textbf{r}_k(X))-g^*(X))\Big]\\
		&\quad+\mathbb{E}\Big[|g^*(\textbf{r}_k(X))-g^*(X)|^2\Big].
	\end{align*}
	We consider the second term of the right hand side of the last equality,
	\begin{align*}
		{}&\ \mathbb{E}\Big[(g_n(\textbf{r}_k(X))-g^*(\textbf{r}_k(X)))(g^*(\textbf{r}_k(X))-g^*(X))\Big]\\
		=&\ \mathbb{E}_{\textbf{r}_k(X)}\Big[\mathbb{E}_{X}\Big[(g_n(\textbf{r}_k(X))-g^*(\textbf{r}_k(X)))(g^*(\textbf{r}_k(X))-g^*(X))\Big|\textbf{r}_k(X)\Big]\Big]\\
		=&\ \mathbb{E}_{\textbf{r}_k(X)}\Big[(g_n(\textbf{r}_k(X))-g^*(\textbf{r}_k(X)))(g^*(\textbf{r}_k(X))-\mathbb{E}[g^*(X)|\textbf{r}_k(X)])\Big]\\
		=&\ 0
	\end{align*}
	where $g^*(\textbf{r}_k(X))=\mathbb{E}[g^*(X)|\textbf{r}_k(X)]$ due to the definition of $g^*(\textbf{r}_k(X))$ and the tower property of conditional expectation. It remains to check that
	$$\mathbb{E}\Big[|g^*(\textbf{r}_k(X))-g^*(X)|^2\Big]\leq \inf_{f\in\mathcal{G}}\mathbb{E}\Big[|f(\textbf{r}_k(X))-g^*(X)|^2\Big].$$
	For any function $f$ s.t $\mathbb{E}\Big[|f(\textbf{r}_k(X))|^2\Big]<+\infty$, one has
	\begin{align*}
		\mathbb{E}\Big[|f(\textbf{r}_k(X))-g^*(X)|^2\Big]&= \mathbb{E}\Big[|f(\textbf{r}_k(X))-g^*(\textbf{r}_k(X))+g^*(\textbf{r}_k(X))-g^*(X)|^2\Big]\\
		&= \mathbb{E}\Big[|f(\textbf{r}_k(X))-g^*(\textbf{r}_k(X))|^2\Big]\\
		&\quad+2\mathbb{E}\Big[(f(\textbf{r}_k(X))-g^*(\textbf{r}_k(X)))(g^*(\textbf{r}_k(X))-g^*(X))\Big]\\
		&\quad+\mathbb{E}\Big[|g^*(\textbf{r}_k(X))-g^*(X)|^2\Big].
	\end{align*}
	Similarly,
	\begin{align*}
		\mathbb{E}\Big[(f(\textbf{r}_k(X))-g^*(\textbf{r}_k(X)))(g^*(\textbf{r}_k(X))-g^*(X))\Big] = 0.
	\end{align*}
	Therefore,
	\begin{align*}
		\mathbb{E}\Big[|f(\textbf{r}_k(X))-g^*(X)|^2\Big]&= \mathbb{E}\Big[|f(\textbf{r}_k(X))-g^*(\textbf{r}_k(X))|^2\Big]+\mathbb{E}\Big[|g^*(\textbf{r}_k(X))-g^*(X)|^2\Big].
	\end{align*}
	As the first term of the right-hand side is nonnegative thus,
	$$\mathbb{E}\Big[|g^*(\textbf{r}_k(X))-g^*(X)|^2\Big]\leq \inf_{f\in\mathcal{G}}\mathbb{E}\Big[|f(\textbf{r}_k(X))-g^*(X)|^2\Big].$$
	Finally, we can conclude that
	\begin{align*}
		\mathbb{E}\Big[|g_n(\textbf{r}_k(X))-g^*(X)|^2\Big]&\leq \mathbb{E}\Big[|g_n(\textbf{r}_k(X))-g^*(\textbf{r}_k(X))|^2\Big]+\inf_{f\in\mathcal{G}}\mathbb{E}\Big[|f(\textbf{r}_k(X))-g^*(X)|^2\Big].
	\end{align*}
\end{proofprop}
We obtain the particular case by restricting $\mathcal{G}$ to be the coordinates of $\textbf{r}_k$, one has
\begin{align*}
	\mathbb{E}\Big[|g_n(\textbf{r}_k(X))-g^*(X)|^2\Big]&\leq \mathbb{E}\Big[|g_n(\textbf{r}_k(X))-g^*(\textbf{r}_k(X))|^2\Big]+\min_{1\leq m\leq M}\mathbb{E}\Big[|r_{k,m}(X)-g^*(X)|^2\Big].
\end{align*}
\QED

\begin{proofprop}
	The procedure of proving this result is indeed the procedure of checking the conditions of Stone's theorem (see, for example, \cite{stone1977} and Chapter 4 of \cite{bookDistributionFree}) which is also used in the classical method by \cite{cobra}. First of all, using the inequality: $(a+b+c)^2\leq 3(a^2+b^2+c^2)$, one has
	\begin{align*}
		\ \mathbb{E}\Big[|g_n(\textbf{r}_k(X))-g^*(\textbf{r}_k(X))|^2\Big]&=\mathbb{E}\Big[\Big|\sum_{i=1}^{\ell}W_{n,i}(X)Y_i-g^*(\textbf{r}_k(X))\Big|^2\Big]\\
		&=\mathbb{E}\Big[\Big|\sum_{i=1}^{\ell}W_{n,i}(X)[Y_i-g^*(\textbf{r}_k(X_i))]\\
		&\quad+\sum_{i=1}^{\ell}W_{n,i}(X)[g^*(\textbf{r}_k(X_i))-g^*(\textbf{r}_k(X))]\\
		&\quad+\sum_{i=1}^{\ell}W_{n,i}(X)g^*(\textbf{r}_k(X))-g^*(\textbf{r}_k(X))\Big|^2\Big]\\
		&\leq 3\mathbb{E}\Big[\Big|\sum_{i=1}^{\ell}W_{n,i}(X)[g^*(\textbf{r}_k(X_i))-g^*(\textbf{r}_k(X))]\Big|^2\Big]\\
		&\quad+3\mathbb{E}\Big[\Big|\sum_{i=1}^{\ell}W_{n,i}(X)[Y_i-g^*(\textbf{r}_k(X_i))]\Big|^2\Big]\\
		&\quad+3\mathbb{E}\Big[\Big|g^*(\textbf{r}_k(X))\sum_{i=1}^{\ell}(W_{n,i}(X)-1)\Big|^2\Big].
	\end{align*}
	The three terms of the right-hand side are denoted by $A.1, A.2$ and $A.3$ respectively, thus one has
	\begin{align*}
		\mathbb{E}\Big[|g_n(\textbf{r}_k(X))-g^*(\textbf{r}_k(X))|^2\Big]&\leq 3(A.1+A.2+A.3).
	\end{align*}
	To prove the result, it is enough to prove that the three terms \textbf{$A.1,A.2$} and \textbf{$A.3$} vanish under the assumptions of \textbf{Proposition~2}. We deal with the first term \textbf{$A.1$} in the following proposition.
	
	\begin{propA}
		Under the assumptions of \textbf{Proposition~2},
		\begin{align*}
			\lim_{\ell\rightarrow+\infty}\mathbb{E}\Big[\Big|\sum_{i=1}^{\ell}W_{n,i}(X)[g^*(\textbf{r}_k(X_i))-g^*(\textbf{r}_k(X))]\Big|^2\Big]=0.
		\end{align*}
	\end{propA}
	\begin{proofA}
		Using Cauchy-Schwarz's inequality, one has
		\begin{align*}
			A.1&=\mathbb{E}\Big[\Big|\sum_{i=1}^{\ell}W_{n,i}(X)[g^*(\textbf{r}_k(X_i))-g^*(\textbf{r}_k(X))]\Big|^2\Big]\\
			&=\mathbb{E}\Big[\Big|\sum_{i=1}^{\ell}\sqrt{W_{n,i}(X)}\sqrt{W_{n,i}(X)}[g^*(\textbf{r}_k(X_i))-g^*(\textbf{r}_k(X))]\Big|^2\Big]\\
			{}&\leq \mathbb{E}\Big[\Big(\sum_{i=1}^{\ell}W_{n,i}(X)\Big)\sum_{i=1}^{\ell}W_{n,i}(X)[g^*(\textbf{r}_k(X_i))-g^*(\textbf{r}_k(X))]^2\Big]\\
			&=\mathbb{E}\Big[\sum_{i=1}^{\ell}W_{n,i}(X)[g^*(\textbf{r}_k(X_i))-g^*(\textbf{r}_k(X))]^2\Big]\\
			&\eqdef A_n.
		\end{align*}
		Note that the regression function $g^*$ satisfies $\mathbb{E}[|g^*(\textbf{r}_k(X))|^2]<+\infty$, thus it can be approximated in $L_2$ sense by a continuous function with compact support named $\tilde{g}$ (see, for example, Theorem A.1 in \cite{bookProbabTheoryPattern}). This means that for any $\varepsilon>0$, there exists a continuous function with compact support $\tilde{g}$ such that,
		$$\mathbb{E}[|g^*(\textbf{r}_k(X))-\tilde{g}(\textbf{r}_k(X))|^2]<\varepsilon.$$
		Thus, one has
		
		\begin{align*}
			A_n&=\mathbb{E}\Big[\sum_{i=1}^{\ell}W_{n,i}(X)[g^*(\textbf{r}_k(X_i))-g^*(\textbf{r}_k(X))]^2\Big]
		\end{align*}
		\begin{align*}
			&\leq 3\mathbb{E}\Big[\sum_{i=1}^{\ell}W_{n,i}(X)[g^*(\textbf{r}_k(X_i))-\tilde{g}(\textbf{r}_k(X_i))]^2\Big]\\
			&\quad+3\mathbb{E}\Big[\sum_{i=1}^{\ell}W_{n,i}(X)[\tilde{g}(\textbf{r}_k(X_i))-\tilde{g}(\textbf{r}_k(X))]^2\Big]\\
			&\quad+3\mathbb{E}\Big[\sum_{i=1}^{\ell}W_{n,i}(X)[\tilde{g}(\textbf{r}_k(X))-g^*(\textbf{r}_k(X))]^2\Big]\\
			&\eqdef 3(A_{n1}+A_{n2}+A_{n3}).
		\end{align*}
		We deal with each term of the last upper bound as follows.
		
		\begin{itemize}
			\item Computation of \textbf{$A_{n3}$}: applying the definition of $\tilde{g}$,
			\begin{align*}
				A_{n3}&=\mathbb{E}\Big[\sum_{i=1}^{\ell}W_{n,i}(X)[\tilde{g}(\textbf{r}_k(X))-g^*(\textbf{r}_k(X))]^2\Big]\\
				&\leq \mathbb{E}\Big[|\tilde{g}(\textbf{r}_k(X))-g^*(\textbf{r}_k(X))|^2\Big]< \varepsilon.
			\end{align*}
			\item Computation of \textbf{$A_{n1}$}: denoted by $\mu$ the distribution of $X$. Thus,
			\begin{align*}
				A_{n1}&=\mathbb{E}\Big[\sum_{i=1}^{\ell}W_{n,i}(X)|g^*(\textbf{r}_k(X_i))-\tilde{g}(\textbf{r}_k(X_i))|^2\Big]\\
				&=\ell\mathbb{E}\Big[W_{n,1}(X)|g^*(\textbf{r}_k(X_1))-\tilde{g}(\textbf{r}_k(X_1))|^2\Big]\\
				&=\ell\mathbb{E}\Big[\frac{K_h(\textbf{r}_k(X)-\textbf{r}_k(X_1))}{\sum_{j=1}^{\ell}K_h(\textbf{r}_k(X)-\textbf{r}_k(X_j))}|g^*(\textbf{r}_k(X_1))-\tilde{g}(\textbf{r}_k(X_1))|^2\Big]\\
				&=\ell\mathbb{E}_{\mathcal{D}_k}\Big[\mathbb{E}_{\{X_j\}_{j=1}^{\ell}}\Big[\int\frac{K_h(\textbf{r}_k(v)-\textbf{r}_k(X_1))}{\sum_{j=1}^{\ell}K_h(\textbf{r}_k(v)-\textbf{r}_k(X_j))}|g^*(\textbf{r}_k(X_1))-\tilde{g}(\textbf{r}_k(X_1))|^2\mu(dv)\Big|\mathcal{D}_k\Big]\Big]\\
				&=\ell\mathbb{E}_{\mathcal{D}_k}\Big[\mathbb{E}_{\{X_j\}_{j=2}^{\ell}}\Big[\int\int|g^*(\textbf{r}_k(u))-\tilde{g}(\textbf{r}_k(u))|^2\times\\
				&\quad\frac{K_h(\textbf{r}_k(v)-\textbf{r}_k(u))}{K_h(\textbf{r}_k(v)-\textbf{r}_k(u))+\sum_{j=2}^{\ell}K_h(\textbf{r}_k(v)-\textbf{r}_k(X_j))}\mu(du)\mu(dv)\Big|\mathcal{D}_k\Big]\Big]\\
				&=\ell\mathbb{E}_{\mathcal{D}_k}\Big[\int|g^*(\textbf{r}_k(u))-\tilde{g}(\textbf{r}_k(u))|^2\times\\
				&\quad\ \mathbb{E}_{\{X_j\}_{j=2}^{\ell}}\Big[\int\frac{K_h(\textbf{r}_k(v)-\textbf{r}_k(u))\mu(dv)}{K_h(\textbf{r}_k(v)-\textbf{r}_k(u))+\sum_{j=2}^{\ell}K_h(\textbf{r}_k(v)-\textbf{r}_k(X_j))}\Big|\mathcal{D}_k\Big]\mu(du)\Big]\\
				&=\ell\mathbb{E}_{\mathcal{D}_k}\Big[\int|g^*(\textbf{r}_k(u))-\tilde{g}(\textbf{r}_k(u))|^2\times I(u,\ell)\mu(du)\Big].
			\end{align*}
			Fubini's theorem (\cite{Folland1999}) is employed to obtain the result of the last bound where the inner conditional expectation is denoted by $I(u,\ell)$. We bound $I(u,\ell)$ using the argument of covering $\mathbb{R}^M$ with a countable family of balls $\mathcal{B}\eqdef\{B_M(x_i,\rho/2):i=1,2,....\}$ and the facts that
			\begin{enumerate}
				\item $\textbf{r}_k(v)\in B_M(\textbf{r}_k(u)+hx_i,h\rho/2)\Rightarrow B_M(\textbf{r}_k(u)+hx_i,h\rho/2)\subset B_M(\textbf{r}_k(v),h\rho)$.
				\item $b\mathds{1}_{\{B_M(0,\rho)\}}(z) < K(z)\leq 1,\forall z\in\mathbb{R}^M.$
			\end{enumerate}
			
			Now, let 
			\begin{itemize}
				\item $A_{i,h}(u)\eqdef\{v\in\mathbb{R}^d:\|\textbf{r}_k(v)-\textbf{r}_k(u)-hx_i\|<h\rho/2\}$.
				\item $B_{i,h}^{\ell}(u)\eqdef\sum_{j=2}^{\ell}\mathds{1}_{\{\|\textbf{r}_k(X_j)-\textbf{r}_k(u)-hx_i\|<h\rho/2\}}$. 
			\end{itemize}
			Thus, one has
			
			\begin{align*}
				I(u,\ell)&\eqdef\mathbb{E}_{\{X_j\}_{j=2}^{\ell}}\Big[\int\frac{K_h(\textbf{r}_k(v)-\textbf{r}_k(u))\mu(dv)}{K_h(\textbf{r}_k(v)-\textbf{r}_k(u))+\sum_{j=2}^{\ell}K_h(\textbf{r}_k(v)-\textbf{r}_k(X_j))}\Big|\mathcal{D}_k\Big]\\
				&\leq\mathbb{E}_{\{X_j\}_{j=2}^{\ell}}\Big[\sum_{i=1}^{+\infty}\int_{v:\|\textbf{r}_k(v)-\textbf{r}_k(u)-hx_i\|<h\rho/2}\\
				&\quad\frac{K_h(\textbf{r}_k(v)-\textbf{r}_k(u))\mu(dv)}{K_h(\textbf{r}_k(v)-\textbf{r}_k(u))+\sum_{j=2}^{\ell}K_h(\textbf{r}_k(v)-\textbf{r}_k(X_j))}\Big|\mathcal{D}_k\Big]\\
				&\leq\mathbb{E}_{\{X_j\}_{j=2}^{\ell}}\Big[\sum_{i=1}^{+\infty}\int_{A_{i,h}(u)}\\
				&\quad\frac{\sup_{z:\|z-hx_i\|<h\rho/2}K_h(z)\mu(dv)}{\sup_{z:\|z-hx_i\|<h\rho/2}K_h(z)+\sum_{j=2}^{\ell}K_h(\textbf{r}_k(v)-\textbf{r}_k(X_j))}\Big|\mathcal{D}_k\Big]\\
				&\leq\mathbb{E}_{\{X_j\}_{j=2}^{\ell}}\Big[\sum_{i=1}^{+\infty}\int_{A_{i,h}(u)}\\
				&\quad\frac{\sup_{z:\|z-hx_i\|<h\rho/2}K_h(z)\mu(dv)}{\sup_{z:\|z-hx_i\|<h\rho/2}K_h(z)+b\sum_{j=2}^{\ell}\mathds{1}_{\{\|\textbf{r}_k(v)-\textbf{r}_k(X_j)\|<h\rho\}}}\Big|\mathcal{D}_k\Big]\\
				&\leq\frac{1}{b}\mathbb{E}_{\{X_j\}_{j=2}^{\ell}}\Big[\sum_{i=1}^{+\infty}\int_{A_{i,h}(u)}\\
				&\quad\frac{\sup_{z:\|z-hx_i\|<h\rho/2}K_h(z)\mu(dv)}{\sup_{z:\|z-hx_i\|<h\rho/2}K_h(z)+\sum_{j=2}^{\ell}\mathds{1}_{\{\|\textbf{r}_k(X_j)-\textbf{r}_k(u)-hx_i\|<h\rho/2\}}}\Big|\mathcal{D}_k\Big]
			\end{align*}
		\begin{align*}
				&\leq\frac{1}{b}\sum_{i=1}^{+\infty}\mathbb{E}_{\{X_j\}_{j=2}^{\ell}}\Big[\frac{\sup_{z:\|z-hx_i\|<h\rho/2}K_h(z)\mu(A_{i,h}(u))}{\sup_{z:\|z-hx_i\|<h\rho/2}K_h(z)+B_{i,h}^{\ell}(u)}\Big|\mathcal{D}_k\Big].
			\end{align*}
			
			Note that $B_{i,h}^{\ell}(u)$ is a binomial random variable $B(\ell-1, \mu(A_{i,h}(u)))$ under the law of $\{X_j\}_{j=2}^{\ell}$. Applying part 1 of lemma~\ref{lem:1}, one has
			
			\begin{align*}
				I(u,\ell)&\leq\frac{1}{b}\sum_{i=1}^{+\infty}\frac{2\sup_{z:\|z-hx_i\|<h\rho/2}K_h(z)\mu(A_{i,h}(u))}{\ell\mu(A_{i,h}(u))}\\
				&\leq\frac{2}{b\ell}\sum_{i=1}^{+\infty}\sup_{w:\|w-x_i\|<\rho/2}K(w)\\
				&=\frac{2}{b\ell}\sum_{i=1}^{+\infty}\sup_{w\in B_M(x_i,\rho/2)}K(w)\\
				&\leq\frac{2}{b\ell}\sum_{i=1}^{+\infty}\sup_{w\in B_M(x_i,\rho/2)}K(w)\\
				&\leq\frac{2}{b\ell\lambda_M(B_M(0,\rho/2))}\sum_{i=1}^{+\infty}\int_{B_M(x_i,\rho/2)}\sup_{w\in B_M(x_i,\rho/2)}K(w)dy\\
				&\leq\frac{2}{b\ell\lambda_M(B_M(0,\rho/2))}\sum_{i=1}^{+\infty}\int_{B_M(x_i,\rho/2)}\sup_{w\in B_M(y,\rho)}K(w)dy\\
				&\leq\frac{2\kappa_M}{b\ell\lambda_M(B_M(0,\rho/2))}\underbrace{\int\sup_{w\in B_M(y,\rho)}K(w)dy}_{=\ \kappa_0\ \text{by \eqref{eq:regular}}}\\
				&\leq\frac{2\kappa_M\kappa_0}{b\ell\lambda_M(B_M(0,\rho))}\eqdef\ \frac{C(b,\rho,\kappa_0,M)}{\ell}<+\infty
			\end{align*}
			where $\lambda_M$ denotes the Lebesque measure on of $\mathbb{R}^M$, $\kappa_M$ denotes the number of balls covering a certain element of $\mathbb{R}^M$, and the constant part is denoted by $C(b,\rho,\kappa_0,M)$ depending on the parameters indicated in the bracket. The last inequality is attained from the fact that the overlapping integrals $\sum_{i=1}^{+\infty}\int_{B_M(x_i,\rho/2)}\sup_{z\in B_M(y,\rho/2)}K(z)dy$ is bounded above by the integral over the entire space $\int\sup_{z\in B_M(y,\rho/2)}K(z)dy$ multiplying by the number of covering balls $k_M$. Therefore,
			
			\begin{align*}
				A_{n1}&\leq\ell\frac{C(b,\rho,\kappa_0,M)}{\ell}\mathbb{E}_{\mathcal{D}_k}\Big[\int|g^*(\textbf{r}_k(u))-\tilde{g}(\textbf{r}_k(u))|^2\mu(du)\Big]\\
				&=C(b,\rho,\kappa_0,M)\mathbb{E}\Big[|\tilde{g}(\textbf{r}_k(X))-g^*(\textbf{r}_k(X))|^2\Big]\\
				&<C(b,\rho,\kappa_0,M)\varepsilon.
			\end{align*}
			
			\item Computation of \textbf{$A_{n2}$}: for any $\delta>0$ one has
			\begin{align*}
				A_{n2}&=\mathbb{E}\Big[\sum_{i=1}^{\ell}W_{n,i}(X)|\tilde{g}(\textbf{r}_k(X_i))-\tilde{g}(\textbf{r}_k(X))|^2\Big]\\
				&=\mathbb{E}\Big[\sum_{i=1}^{\ell}W_{n,i}(X)|\tilde{g}(\textbf{r}_k(X_i))-\tilde{g}(\textbf{r}_k(X))|^2\mathds{1}_{\{\|\textbf{r}_{k}(X_i)-\textbf{r}_k(X)\|\geq\delta\}}\Big]\\
				&\quad+\mathbb{E}\Big[\sum_{i=1}^{\ell}W_{n,i}(X)|\tilde{g}(\textbf{r}_k(X_i))-\tilde{g}(\textbf{r}_k(X))|^2\mathds{1}_{\{\|\textbf{r}_{k}(X_i)-\textbf{r}_k(X)\|<\delta\}}\Big]\\
				&\leq4\sup_{u\in\mathbb{R}^d}|\tilde{g}(\textbf{r}_k(u))|^2\mathbb{E}\Big[\sum_{i=1}^{\ell}W_{n,i}(X)\mathds{1}_{\{\|\textbf{r}_{k}(X_i)-\textbf{r}_k(X)\|\geq\delta\}}\Big]\\
				&\quad+ \sup_{u,v\in\mathbb{R}^d:\|\textbf{r}_{k}(u)-\textbf{r}_{k}(v)\|<\delta}|\tilde{g}(\textbf{r}_k(u))-\tilde{g}(\textbf{r}_k(v))|^2
			\end{align*}
			Using the uniform continuity of $\tilde{g}$, the second term of the upper bound of $A_{n2}$ tends to $0$ when $\delta$ tends $0$. Thus, we only need to prove that the first term of this upper bound also tends to $0$. We follow a similar procedure as in the previous part:
			
			\begin{align*}
				{}&\mathbb{E}\Big[\sum_{i=1}^{\ell}W_{n,i}(X)\mathds{1}_{\{\|\textbf{r}_{k}(X_i)-\textbf{r}_k(X)\|\geq\delta\}}\Big]\\
				=&\ \mathbb{E}_{\mathcal{D}_k}\Big[\sum_{i=1}^{\ell}\mathbb{E}_{X,\{X_j\}_{j=1}^{\ell}}\Big[W_{n,i}(X)\mathds{1}_{\{\|\textbf{r}_{k}(X)-\textbf{r}_k(X_i)\|\geq\delta\}}\Big|\mathcal{D}_k\Big]\Big]\\
				=&\ \mathbb{E}_{\mathcal{D}_k}\Big[\sum_{i=1}^{\ell}\mathbb{E}_{\{X_j\}_{j=1}^{\ell}}\Big[\int\frac{K_h(\textbf{r}_k(v)-\textbf{r}_k(X_i))\mathds{1}_{\{\|\textbf{r}_{k}(v)-\textbf{r}_k(X_i)\|\geq\delta\}}}{\sum_{j=1}^{\ell}K_h(\textbf{r}_k(v)-\textbf{r}_k(X_j))}\mu(dv)\Big|\mathcal{D}_k\Big]\Big]\\
				=&\ \ell\mathbb{E}_{\mathcal{D}_k}\Big[\mathbb{E}_{\{X_j\}_{j=2}^{\ell}}\Big[\int\int\frac{K_h(\textbf{r}_k(v)-\textbf{r}_k(u))\mathds{1}_{\{\|\textbf{r}_{k}(v)-\textbf{r}_k(u)\|\geq\delta\}}\mu(du)\mu(dv)}{K_h(\textbf{r}_k(v)-\textbf{r}_k(u))+\sum_{j=2}^{\ell}K_h(\textbf{r}_k(v)-\textbf{r}_k(X_j))}\Big|\mathcal{D}_k\Big]\Big]\\
				=&\ \ell\mathbb{E}_{\mathcal{D}_k}\Big[\int J(u,\ell)\mu(du)\Big].
			\end{align*}
			Fubini's theorem is applied to obtain the last equation where for any $u\in\mathbb{R}^d$,
			\begin{align*}
				J(u,\ell)&\eqdef\mathbb{E}_{\{X_j\}_{j=2}^{\ell}}\Big[\int\frac{K_h(\textbf{r}_k(v)-\textbf{r}_k(u))\mathds{1}_{\{\|\textbf{r}_{k}(v)-\textbf{r}_k(u)\|\geq\delta\}}\mu(dv)}{K_h(\textbf{r}_k(v)-\textbf{r}_k(u))+\sum_{j=2}^{\ell}K_h(\textbf{r}_k(v)-\textbf{r}_k(X_j))}\Big|\mathcal{D}_k\Big]\\
				&\leq\mathbb{E}_{\{X_j\}_{j=2}^{\ell}}\Big[\sum_{i=1}^{+\infty}\int_{v:\|\textbf{r}_k(v)-\textbf{r}_k(u)-hx_i\|<h\rho/2}\\
				&\quad\ \frac{K_h(\textbf{r}_k(v)-\textbf{r}_k(u))\mathds{1}_{\{\|\textbf{r}_{k}(v)-\textbf{r}_k(u)\|\geq\delta\}}}{K_h(\textbf{r}_k(v)-\textbf{r}_k(u))+\sum_{j=2}^{\ell}K_h(\textbf{r}_k(v)-\textbf{r}_k(X_j))}\mu(dv)\Big|\mathcal{D}_k\Big]\\
				&\leq\mathbb{E}_{\{X_j\}_{j=2}^{\ell}}\Big[\sum_{i=1}^{+\infty}\int_{A_{i,h}(u)}\\
				&\quad\ \frac{\sup_{z:\|z-hx_i\|<h\rho/2}K_h(z)\mathds{1}_{\{\|z\|\geq\delta\}}}{\sup_{z:\|z-hx_i\|<h\rho/2}K_h(z)+\sum_{j=2}^{\ell}K_h(\textbf{r}_k(v)-\textbf{r}_k(X_j))}\mu(dv)\Big|\mathcal{D}_k\Big]\\
				&\leq\sum_{i=1}^{+\infty}\sup_{z:\|z-hx_i\|<h\rho/2}K_h(z)\mathds{1}_{\{\|z\|\geq\delta\}}\times\mathbb{E}_{\{X_j\}_{j=2}^{\ell}}\Big[\int_{A_{i,h}(u)}\\
				&\quad\ \frac{\mu(dv)}{\sup_{z:\|z-hx_i\|<h\rho/2}K_h(z)+b\sum_{j=2}^{\ell}\mathds{1}_{\{\|\textbf{r}_k(X_j)-\textbf{r}_k(v)\|<h\rho\}}}\Big|\mathcal{D}_k\Big]\\
				&\leq\sum_{i=1}^{+\infty}\sup_{z:\|z-hx_i\|<h\rho/2}K_h(z)\mathds{1}_{\{\|z\|\geq\delta\}}\times\mathbb{E}_{\{X_j\}_{j=2}^{\ell}}\Big[\int_{A_{i,h}(u)}\\
				&\quad\ \mathbb{E}_{\{X_j\}_{j=2}^{\ell}}\Big[\frac{\mu(dv)}{\sup_{z:\|z-hx_i\|<h\rho/2}K_h(z)+b\sum_{j=2}^{\ell}\mathds{1}_{\{\|\textbf{r}_k(X_j)-\textbf{r}_k(u)-hx_i\|<h\rho/2\}}}\Big|\mathcal{D}_k\Big]\\
				&\leq\sum_{i=1}^{+\infty}\sup_{z:\|z-hx_i\|<h\rho/2}K_h(z)\mathds{1}_{\{\|z\|\geq\delta\}}\mu(A_{i,h}(u))\times\\
				&\quad\ \frac{1}{b}\mathbb{E}_{\{X_j\}_{j=2}^{\ell}}\Big[\frac{1}{\sup_{z:\|z-hx_i\|<h\rho/2}K_h(z)+B_{i,h}^{\ell}(u)}\Big|\mathcal{D}_k\Big]\\
				&\leq\frac{1}{b}\sum_{i=1}^{+\infty}\frac{2\sup_{z:\|z-hx_i\|<h\rho/2}K_h(z)\mu(A_{i,h}(u))\mathds{1}_{\{\|z\|\geq\delta\}}}{\ell \mu(A_{i,h}(u))}\\
				&\leq\frac{2}{b\ell}\sum_{i=1}^{+\infty}\sup_{w:\|w-x_i\|<\rho/2}K(w)\mathds{1}_{\{\|w\|\geq\delta/h\}}.
			\end{align*}
			Thus, one has
			\begin{align*}
				\mathbb{E}\Big[\sum_{i=1}^{\ell}W_{n,i}(X)\mathds{1}_{\{\|\textbf{r}_{k}(X_i)-\textbf{r}_k(X)\|\geq\delta\}}\Big]
				&\leq \ell\frac{2}{b\ell}\sum_{i=1}^{+\infty}\sup_{w\in B_M(x_i,\rho/2)}K(w)\mathds{1}_{\{\|w\|\geq\delta/h\}}
			\end{align*}
			When both $h\to0$ and $\delta\to0$ satisfying $\delta/h\to+\infty$, the upper bound series converges to zero. Indeed, it is a non-negative convergent series thanks to the proof of $I(u,l)$ in the previous part. Moreover, the general term of the series, $s_k=\sup_{w\in B_M(x_k,\rho/2)}K(w)\mathds{1}_{\{\|w\|\geq\delta/h\}}$, satisfying $\lim_{\delta/h\to+\infty} s_k=0$ for all $k\geq 1$. Therefore, this series converges to zero when $h\to0,\delta\to0$ such that $\delta/h\to+\infty$.
		\end{itemize}
		In conclusion, when $\ell\to+\infty$ and $\varepsilon,h,\delta\to0$ such that $\delta/h\to+\infty$, all the three terms of the upper bound of $A_n$ tend to $0$, so does $A_n$.
		
		\QED
	\end{proofA}
	
	\begin{propA}
		Under the assumptions of \textbf{Proposition~2},
		\begin{align*}
			\lim_{\ell\rightarrow+\infty}\mathbb{E}\Big[\Big|\sum_{i=1}^{\ell}W_{n,i}(X)[Y_i-g_n(\textbf{r}_k(X_i))]\Big|^2\Big]=0.
		\end{align*}
	\end{propA}
	
	\begin{proofA} Using the independence between $(X_i,Y_i)$ and $(X_j,Y_j)$ for all $i\neq j$, one has 
		\begin{align*}
			A.2&=\mathbb{E}\Big[\Big|\sum_{i=1}^{\ell}W_{n,i}(X)[Y_i-g_n(\textbf{r}_k(X_i))]\Big|^2\Big]\\
			&=\sum_{1\leq i,j\leq\ell}\mathbb{E}\Big[W_{n,i}(X)W_{n,j}(X)[Y_i-g_n(\textbf{r}_k(X_i))][Y_j-g_n(\textbf{r}_k(X_j))]\Big]\\
			&=\ \mathbb{E}\Big[\sum_{i=1}^{\ell}W_{n,i}^2(X)|Y_i-g_n(\textbf{r}_k(X_i))|^2\Big]=\mathbb{E}\Big[\sum_{i=1}^{\ell}W_{n,i}^2(X)\sigma^2(\textbf{r}_k(X_i))\Big]
		\end{align*}
		where 
		$$\sigma^2(\textbf{r}_k(x))\eqdef\mathbb{E}[(Y_i-g_n(\textbf{r}_k(X_i)))^2|\textbf{r}_k(x)].$$
		Thus, based on the assumption of $X$ and $Y$ we have $\sigma^2\in L_1(\mu)$. Therefore, $\sigma^2$ can be approximated in $L_1$ sense i.e., for any $\varepsilon>0, \exists\tilde{\sigma}^2$ a continuous function with compact support such that
		$$\mathbb{E}[|\sigma^2(\textbf{r}_k(X))-\tilde{\sigma}^2(\textbf{r}_k(X))|]<\varepsilon.$$
		Thus, one has
		\begin{align*}
			A.2&\leq\mathbb{E}\Big[\sum_{i=1}^{\ell}W_{n,i}^2(X)\tilde{\sigma}^2(\textbf{r}_k(X_i))\Big]+\mathbb{E}\Big[\sum_{i=1}^{\ell}W_{n,i}^2(X)|\sigma^2(\textbf{r}_k(X_i))-\tilde{\sigma}^2(\textbf{r}_k(X_i))|\Big]\\
			&\leq\sup_{u\in\mathbb{R}^d}|\tilde{\sigma}^2(\textbf{r}_k(u))|\mathbb{E}\Big[\sum_{i=1}^{\ell}W_{n,i}^2(X)\Big]+\mathbb{E}\Big[\sum_{i=1}^{\ell}W_{n,i}^2(X)|\sigma^2(\textbf{r}_k(X_i))-\tilde{\sigma}^2(\textbf{r}_k(X_i))|\Big].
		\end{align*}
		Using similar argument as in the case of $A_{n1}$ and the fact that $W_{n,i}(x)\leq 1,\forall i=1,2,...,\ell$, thus for any $\varepsilon>0$, one has
		\begin{align*}
			\mathbb{E}\Big[\sum_{i=1}^{\ell}W_{n,i}^2(X)|\sigma^2(\textbf{r}_k(X_i))-\tilde{\sigma}^2(\textbf{r}_k(X_i))|\Big]&\leq\mathbb{E}\Big[\sum_{i=1}^{\ell}W_{n,i}(X)|\sigma^2(\textbf{r}_k(X_i))-\tilde{\sigma}^2(\textbf{r}_k(X_i))|\Big]\\
			&<C(b,\rho,\kappa_0,M)\varepsilon.
		\end{align*}
		Therefore, it remains to prove that $\mathbb{E}[\sum_{i=1}^{\ell}W_{n,i}^2(X)]\to0$ as $\ell\to+\infty$. As $b\mathds{1}_{\{B_M(0,\rho)\}}(z) < K(z)\leq 1,\forall z\in\mathbb{R}^M$ with the convention of $0/0=0$, for a fixed $\delta>0$, one has
		\begin{align}
			\sum_{i=1}^{\ell}W_{n,i}^2(X)&=\sum_{i=1}^{\ell}\Big(\frac{K_h(\textbf{r}_k(X)-\textbf{r}_k(X_i))}{\sum_{j=1}^{\ell}K_h(\textbf{r}_k(X)-\textbf{r}_k(X_j))}\Big)^2\nonumber\\
			&\leq\frac{\sum_{i=1}^{\ell}K_h(\textbf{r}_k(X)-\textbf{r}_k(X_i))}{\Big(\sum_{j=1}^{\ell}K_h(\textbf{r}_k(X)-\textbf{r}_k(X_j))\Big)^2}\nonumber\\
			&\leq\min\Big\{\delta,\frac{\mathds{1}_{\{\sum_{j=1}^{\ell}K_h(\textbf{r}_k(X)-\textbf{r}_k(X_j))>0\}}}{\sum_{j=1}^{\ell}K_h(\textbf{r}_k(X)-\textbf{r}_k(X_j))}\Big\}\nonumber\\
			&\leq\min\Big\{\delta,\frac{\mathds{1}_{\{\sum_{j=1}^{\ell}\mathds{1}_{\{\|\textbf{r}_k(X)-\textbf{r}_k(X_j)\|<h\rho\}}>0\}}}{b\sum_{j=1}^{\ell}\mathds{1}_{\{\|\textbf{r}_k(X)-\textbf{r}_k(X_j)\|<h\rho\}}}\nonumber\\
			&\leq\delta+\frac{\mathds{1}_{\{\sum_{j=1}^{\ell}\mathds{1}_{\{\|\textbf{r}_k(X)-\textbf{r}_k(X_j)\|<h\rho\}}>0\}}}{b\sum_{j=1}^{\ell}\mathds{1}_{\{\|\textbf{r}_k(X)-\textbf{r}_k(X_j)\|<h\rho\}}}.\label{eq:boundW2}
		\end{align}
		Therefore, it is enough to show that
		$$\mathbb{E}\Big[\frac{\mathds{1}_{\{\sum_{j=1}^{\ell}\mathds{1}_{\{\|\textbf{r}_k(X)-\textbf{r}_k(X_j)\|<h\rho\}}>0\}}}{\sum_{j=1}^{\ell}\mathds{1}_{\{\|\textbf{r}_k(X)-\textbf{r}_k(X_j)\|<h\rho\}}}\Big]\xrightarrow{\ell\to+\infty}0.$$
		One has
		\begin{align*}
			{}&\ \ \ \ \mathbb{E}\Big[\frac{\mathds{1}_{\{\sum_{j=1}^{\ell}\mathds{1}_{\{\|\textbf{r}_k(X)-\textbf{r}_k(X_j)\|<h\rho\}}>0\}}}{\sum_{j=1}^{\ell}\mathds{1}_{\{\|\textbf{r}_k(X)-\textbf{r}_k(X_j)\|<h\rho\}}}\Big]\\
			&\leq\mathbb{E}\Big[\frac{\mathds{1}_{\{\sum_{j=1}^{\ell}\mathds{1}_{\{\|\textbf{r}_k(X)-\textbf{r}_k(X_j)\|<h\rho\}}>0\}}}{\sum_{j=1}^{\ell}\mathds{1}_{\{\|\textbf{r}_k(X)-\textbf{r}_k(X_j)\|<h\rho\}}}\mathds{1}_{\{\textbf{r}_k(X)\in B\}}\Big]+\mu(\{v\in\mathbb{R}^d:\textbf{r}_k(v)\in B^c\})\\
			&=\mathbb{E}\Big[\mathds{1}_{\{\textbf{r}_k(X)\in B\}}\mathbb{E}\Big[\frac{\mathds{1}_{\{\sum_{j=1}^{\ell}\mathds{1}_{\{\|\textbf{r}_k(X)-\textbf{r}_k(X_j)\|<h\rho\}}>0\}}}{\sum_{j=1}^{\ell}\mathds{1}_{\{\|\textbf{r}_k(X)-\textbf{r}_k(X_j)\|<h\rho\}}}\Big|X\Big]\Big]+\mu(\{v\in\mathbb{R}^d:\textbf{r}_k(v)\in B^c\})\\
			&\leq2\mathbb{E}\Big[\frac{\mathds{1}_{\{\textbf{r}_k(X)\in B\}}}{(\ell+1)\mu(\{v\in\mathbb{R}^d:\|\textbf{r}_k(v)-\textbf{r}_k(X)\|<h\rho\})}\Big]+\mu(\{v\in\mathbb{R}^d:\textbf{r}_k(v)\in B^c\})
		\end{align*}
		where $B$ is a $M$-dimensional ball centered at the origin chosen so that the second term $\mu(\{v\in\mathbb{R}^d:\textbf{r}_k(v)\in B^c\})$ is small. The last inequality is attained by applying part 2 of lemma 1. Moreover, as $\textbf{r}_k=(\textbf{r}_{k,m})_{m=1}^M$ is bounded then there exists a finite number of balls in $\mathcal{B}=\{B_M(x_j,h\rho/2):j=1,2,...\}$ such that $B$ is contained in the union of these balls i.e., $\exists I_{h,M}$ finite, such that $B\subset\cup_{j\in I_{h,M}}B_M(x_j,h\rho/2)$. 
		
		\begin{align}
			{}&\quad\ \mathbb{E}\Big[\frac{\mathds{1}_{\{\textbf{r}_k(X)\in B\}}}{(\ell+1)\mu(\{v\in\mathbb{R}^d:\|\textbf{r}_k(v)-\textbf{r}_k(X)\|<h\rho\})}\Big]\nonumber\\
			&\leq\sum_{j\in I_{h,M}}\int_{u:\|\textbf{r}_k(u)-x_j\|<h\rho/2}\frac{\mu(du)}{(\ell+1)\mu(\{v\in\mathbb{R}^d:\|\textbf{r}_k(v)-\textbf{r}_k(u)\|<h\rho\})}\nonumber\\
			&\quad+\mu(\{v\in\mathbb{R}^d:\textbf{r}_k(v)\in B^c\})\nonumber\\
			&\leq\sum_{j\in I_{h,M}}\int_{u:\|\textbf{r}_k(u)-x_j\|<h\rho/2}\frac{\mu(du)}{(\ell+1)\mu(\{v\in\mathbb{R}^d:\|\textbf{r}_k(v)-x_j\|<h\rho/2\})}\nonumber\\
			&\quad+\mu(\{v\in\mathbb{R}^d:\textbf{r}_k(v)\in B^c\})\nonumber\\
			&=\sum_{j\in I_{h,M}}\frac{\mu(\{u\in\mathbb{R}^d:\|\textbf{r}_k(u)-x_j\|<h\rho/2\})}{(\ell+1)\mu(\{v\in\mathbb{R}^d:\|\textbf{r}_k(v)-x_j\|<h\rho/2\})}+\mu(\{v\in\mathbb{R}^d:\textbf{r}_k(v)\in B^c\})\nonumber\\
			&=\frac{|I_{h,M}|}{\ell+1}+\mu(\{v\in\mathbb{R}^d:\textbf{r}_k(v)\in B^c\})\nonumber\\
			&\leq\frac{C_0}{h^M(\ell+1)}+\mu(\{v\in\mathbb{R}^d:\textbf{r}_k(v)\in B^c\})\label{eq:bound1}\\
			&\xrightarrow[h^M\ell\to+\infty]{\ell\to+\infty,h\to0}\mu(\{v\in\mathbb{R}^d:\textbf{r}_k(v)\in B^c\}). \nonumber
		\end{align}
		It is easy to check the following fact,
		\begin{equation}
			\label{eq:boundJhM}
			|I_{h,M}|\leq\frac{C_0}{h^M}\ \text{for some }C_0>0.
		\end{equation}
		To prove inequality~\eqref{eq:boundJhM}, we consider again the cover $\mathcal{B}=\{B_M(x_j,h\rho/2):j=1,2,...\}$ of $\mathbb{R}^M$. For any $\rho>0$ fixed and $h>0$, note that
		the covering number $|I_{h,M}|$ is proportional to the ratio between the volume of $B$ and the volume of the ball $B_M(0,h\rho/2)$ i.e.,
		\begin{align*}
			|I_{h,M}|&\propto \frac{\text{Vol}(B)}{\text{Vol}(B_M(0,h\rho/2))}\\
			&\propto \frac{\text{Vol}(B)}{(h\rho/2)^M}\\
			&\leq  \frac{C_0}{h^M}
		\end{align*}
		for some positive constant $C_0$ proportional to the volume of $B$. Finally, we can conclude the proof of the proposition as we can choose $B$ such that $\mu(\{v\in\mathbb{R}^d:\textbf{r}_k(v)\in B^c\})=0$ using the boundedness of the basic regressors.
		\begin{remark}
			The assumption on the boundedness of the constructed estimators is crucial. This assumption allows us to choose a ball B which can be covered using a finite number $|I_{h,M}|$ of balls $B_M(x_j,h\rho/2)$, therefore makes it possible to prove the result of this proposition for this class of regular kernels. Note that for the class of compactly supported kernels, it is easy to obtain such a result directly from the begging of the evaluation of each integral (see, for example, Chapter 5 of \cite{bookDistributionFree}).
		\end{remark}
		
		\QED
	\end{proofA}
	\begin{propA}
		Under the assumptions of \textbf{Proposition~2},
		\begin{align*}
			\lim_{\ell\rightarrow+\infty}\mathbb{E}\Big[\Big|g^*(\textbf{r}_k(X))\Big(\sum_{i=1}^{\ell}W_{n,i}(X)-1\Big)\Big|^2\Big]=0.
		\end{align*}
	\end{propA}
	\begin{proofA}
		Note that $|\sum_{i=1}^{\ell}W_{n,i}(X)-1|\leq 1$ thus one has
		$$\Big|g^*(\textbf{r}_k(X))\Big(\sum_{i=1}^{\ell}W_{n,i}(X)-1\Big)\Big|^2\leq|g^*(\textbf{r}_k(X))|^2.$$
		Consequently, by Lebesque's dominated convergence theorem, to prove this proposition, it is enough to show that $\sum_{i=1}^{\ell}W_{n,i}(X)\to1$ almost surely.
		Note that $1-\sum_{i=1}^{\ell}W_{n,i}(X)=\mathds{1}_{\{\sum_{i=1}^{\ell}K_h(\textbf{r}_k(X)-\textbf{r}_k(X_i))=0\}}$ therefore,
		\begin{align*}
			\mathbb{P}\Big[\sum_{i=1}^{\ell}W_{n,i}(X)\neq 1\Big]&=\mathbb{P}\Big[\sum_{i=1}^{\ell}K_h(\textbf{r}_k(X)-\textbf{r}_k(X_i))=0\Big]\\
			&\leq\mathbb{P}\Big(\sum_{j=1}^{\ell}\mathds{1}_{\{\|\textbf{r}_k(X)-\textbf{r}_k(X_j)\|<h\rho\}}=0\Big)\\
			&= \int\mathbb{P}\Big(\sum_{j=1}^{\ell}\mathds{1}_{\{\|\textbf{r}_k(x)-\textbf{r}_k(X_j)\|<h\rho\}}=0\Big)\mu(dx)\\
			&= \int\mathbb{P}\Big(\cap_{j=1}^{\ell}\{\|\textbf{r}_k(x)-\textbf{r}_k(X_j)\|\geq h\rho\}\Big)\mu(dx)\\
			&= \int\Big[1-\mathbb{P}\Big(\{\|\textbf{r}_k(x)-\textbf{r}_k(X_1)\|< h\rho\}\Big)\Big]^{\ell}\mu(dx)
		\end{align*}
		\begin{align*}
			&= \int\Big[1-\mu\Big(\{v\in\mathbb{R}^d:\|\textbf{r}_k(x)-\textbf{r}_k(v)\|< h\rho\}\Big)\Big]^{\ell}\mu(dx)\\
			&\leq \int e^{-\ell\mu(A_h(x))}\mu(dx)\\
			&=\int e^{-\ell\mu(A_h(x))}\mathds{1}_{\{\textbf{r}_k(x)\in B\}}\mu(dx)+\mu(\{v\in\mathbb{R}^d:\textbf{r}_k(v)\in B^c\})\\
			&\leq\frac{\max_{u}\{ue^{-u}\}}{\ell}\int \frac{\mathds{1}_{\{\textbf{r}_k(x)\in B\}}}{\mu(A_h(x))}\mu(dx)+\mu(\{v\in\mathbb{R}^d:\textbf{r}_k(v)\in B^c\})
		\end{align*}
		where
		\begin{equation}
			\label{eq:Axh}
			A_h(x)\eqdef\{v\in\mathbb{R}^d:\|\textbf{r}_k(x)-\textbf{r}_k(v)\|< h\rho\}.
		\end{equation}
		Therefore,
		\begin{align*}
			\mathbb{P}\Big[\sum_{i=1}^{\ell}W_{n,i}(X)\neq 1\Big]&\leq\frac{e^{-1}}{\ell}\mathbb{E}\Big[\frac{\mathds{1}_{\{\textbf{r}_k(X)\in B\}}}{\mu(\{v\in\mathbb{R}^d:\|\textbf{r}_k(v)-\textbf{r}_k(X)\|<h\rho\})}\Big]\\
			&\quad+\mu(\{v\in\mathbb{R}^d:\textbf{r}_k(v)\in B^c\}).
		\end{align*}
		Following the same procedure as in the proof of $A.2$ we obtain the desire result.
	\end{proofA}
	\QED
\end{proofprop}

\begin{proofthm}
	Choose a new observation $x\in\mathbb{R}^d$, given the training data $\mathcal{D}_k$ and the predictions $\{\textbf{r}_k(X_p)\}_{p=1}^{\ell}$ on $\mathcal{D}_{\ell}$, taking expectation with respect to the response variables $\{Y_p^{(\ell)}\}_{p=1}^{\ell}$, it is easy to check that
	\begin{align*}
		&{}\ \ \ \ \mathbb{E}[|g_n(\textbf{r}_k(x))-g^*(\textbf{r}_k(x))|^2|\{\textbf{r}_k(X_p)\}_{p=1}^{\ell},\mathcal{D}_k]\\
		&=\mathbb{E}\Big[\Big|g_n(\textbf{r}_k(x))-\mathbb{E}[g_n(\textbf{r}_k(x))|\{\textbf{r}_k(X_p)\}_{p=1}^{\ell},\mathcal{D}_k]\\
		&\quad+\mathbb{E}[g_n(\textbf{r}_k(x))|\{\textbf{r}_k(X_p)\}_{p=1}^{\ell},\mathcal{D}_k]-g^*(\textbf{r}_k(x))\Big|^2\Big|\{\textbf{r}_k(X_p)\}_{p=1}^{\ell},\mathcal{D}_k\Big]\\
		&=\mathbb{E}[|g_n(\textbf{r}_k(x))-\mathbb{E}[g_n(\textbf{r}_k(x))|\{\textbf{r}_k(X_p)\}_{p=1}^{\ell},\mathcal{D}_k]|^2|\{\textbf{r}_k(X_p)\}_{p=1}^{\ell},\mathcal{D}_k]\\
		&\quad+|g^*(\textbf{r}_k(x))-\mathbb{E}[g_n(\textbf{r}_k(x))|\{\textbf{r}_k(X_p)\}_{p=1}^{\ell},\mathcal{D}_k]|^2\\
		&\eqdef E_1+E_2.
	\end{align*}
	On one hand by using the independence between $Y_i$ and $(Y_j,X_j)$ for all $i\neq j$, we develop the square and obtain for any $\delta>0$:
	
	\begin{align*}
		E_1&\eqdef\mathbb{E}\Big[\Big|g_n(\textbf{r}_k(x))-\mathbb{E}[g_n(\textbf{r}_k(x))|\{\textbf{r}_k(X_p)\}_{p=1}^{\ell},\mathcal{D}_k]\Big|^2\Big|\{\textbf{r}_k(X_p)\}_{p=1}^{\ell},\mathcal{D}_k\Big]\\
		&=\mathbb{E}\Big[\Big|\sum_{i=1}^{\ell}W_{n,i}(x)(Y_i-\mathbb{E}[Y_i|\textbf{r}_k(X_i)])\Big|^2\Big|\{\textbf{r}_k(X_p)\}_{p=1}^{\ell},\mathcal{D}_k\Big]\\
		&=\mathbb{E}\Big[\sum_{i=1}^{\ell}W_{n,i}^2(x)(Y_i-\mathbb{E}[Y_i|\textbf{r}_k(X_i)])^2\Big|\{\textbf{r}_k(X_p)\}_{p=1}^{\ell},\mathcal{D}_k\Big]\\
		&=\sum_{i=1}^{\ell}W_{n,i}^2(x)\mathbb{E}_{Y_i}[(Y_i-\mathbb{E}[Y_i|\textbf{r}_k(X_i)])^2|\textbf{r}_k(X_i)]\\
		&=\mathbb{V}[Y_1|\textbf{r}_k(X_1)]\sum_{i=1}^{\ell}W_{n,i}^2(x)\\
		&\overset{(\ref{eq:boundW2})}{\leq} \frac{4R^2}{b}\Big(\delta+\frac{\mathds{1}_{\{\sum_{j=1}^{\ell}\mathds{1}_{\{\|\textbf{r}_k(x)-\textbf{r}_k(X_j)\|<h\rho\}}>0\}}}{\sum_{j=1}^{\ell}\mathds{1}_{\{\|\textbf{r}_k(x)-\textbf{r}_k(X_j)\|<h\rho\}}}\Big)
	\end{align*}
	where the notation $\mathbb{V}(Z)$ stands for the variance of a random variable $Z$. Therefore, using the result of inequality~(\ref{eq:bound1}), one has
	\begin{equation}
		\label{eq:boundE1}
		\mathbb{E}(E_1)\leq \frac{4R^2}{b}\Big(\delta+\frac{C_0}{h^M(\ell+1)}\Big)
	\end{equation}
	for some $C_0>0$. On the other hand, set 
	\begin{itemize}
		\item[--] $C_{h}^{\ell}(x)\eqdef\sum_{j=1}^{\ell}\mathds{1}_{\{\|\textbf{r}_k(X_j)-\textbf{r}_k(x)\|<h\rho\}}.$
		\item[--] $D_{h}^{\ell}(x)\eqdef\sum_{j=1}^{\ell}K_h({r}_k(X_j)-\textbf{r}_k(x)).$
	\end{itemize}
	The second term $E_2$ is much harder to control as it depends on $g^*(\textbf{r}_k(.))$, that is why a weak smoothness assumption of the theorem is made. Using this assumption and Jensen's inequality (\cite{jensen1906fonctions}), one has
	\begin{align*}
		E_2&\eqdef\Big|g^*(\textbf{r}_k(x))-\mathbb{E}[g_n(\textbf{r}_k(x))|\{\textbf{r}_k(X_p)\}_{p=1}^{\ell},\mathcal{D}_k]\Big|^2\\
		&=\Big(\sum_{i=1}^{\ell}W_{n,i}(X)(g^*(\textbf{r}_k(x))-\mathbb{E}[Y_i|\textbf{r}_k(X_i)])\Big)^2\mathds{1}_{\{D_h^{\ell}(x)>0\}}+(g^*(\textbf{r}_k(x)))^2\mathds{1}_{\{D_h^{\ell}(x)=0\}}\\
		&\overset{(\text{Jensen})}{\leq\ \ \ \ \ }\sum_{i=1}^{\ell}W_{n,i}(x)(g^*(\textbf{r}_k(x))-\mathbb{E}[Y_i|\textbf{r}_k(X_i)])^2\mathds{1}_{\{D_h^{\ell}(x)>0\}}+(g^*(\textbf{r}_k(x)))^2\mathds{1}_{\{D_h^{\ell}(x)=0\}}
	\end{align*}
	\begin{align*}
		&\leq\sum_{i=1}^{\ell}\frac{K_h(\textbf{r}_k(x)-\textbf{r}_k(X_i))(g^*(\textbf{r}_k(x))-g^*(\textbf{r}_k(X_i)))^2}{\sum_{j=1}^{\ell}K_h(\textbf{r}_k(x)-\textbf{r}_k(X_j))}\mathds{1}_{\{D_h^{\ell}(x)>0\}}+(g^*(\textbf{r}_k(x)))^2\mathds{1}_{\{D_h^{\ell}(x)=0\}}\\
		&\leq L^2\sum_{i=1}^{\ell}\frac{K_h(\textbf{r}_k(x)-\textbf{r}_k(X_i))\|\textbf{r}_k(x)-\textbf{r}_k(X_i)\|^2}{\sum_{j=1}^{\ell}K_h(\textbf{r}_k(x)-\textbf{r}_k(X_j))}\mathds{1}_{\{D_h^{\ell}(x)>0\}}+(g^*(\textbf{r}_k(x)))^2\mathds{1}_{\{D_h^{\ell}(x)=0\}}\\
		&\leq L^2\Big[\sum_{i=1}^{\ell}\frac{K_h(\textbf{r}_k(x)-\textbf{r}_k(X_i))\|\textbf{r}_k(x)-\textbf{r}_k(X_i)\|^2\mathds{1}_{\{\|\textbf{r}_k(x)-\textbf{r}_k(X_i)\|<R_Kh^{\beta}\}}}{\sum_{j=1}^{\ell}K_h(\textbf{r}_k(x)-\textbf{r}_k(X_j))}\\
		&\quad+\sum_{i=1}^{\ell}\frac{K_h(\textbf{r}_k(x)-\textbf{r}_k(X_i))\|\textbf{r}_k(x)-\textbf{r}_k(X_i)\|^2\mathds{1}_{\{\|\textbf{r}_k(x)-\textbf{r}_k(X_i)\|\geq R_Kh^{\beta}\}}}{\sum_{j=1}^{\ell}K_h(\textbf{r}_k(x)-\textbf{r}_k(X_j))}\Big]\mathds{1}_{\{D_h^{\ell}(x)>0\}}\\
		&\quad+(g^*(\textbf{r}_k(x)))^2\mathds{1}_{\{C_h^{\ell}(x)=0\}}\\
		&\eqdef E_2^1+E_2^2+E_2^3.
	\end{align*}
	for any $\beta>0$ chosen arbitrarily at this point. Now, we bound the expectation of the three terms of the last inequality. 
	
	\begin{itemize}
		\item Firstly, $E_2^1$ can be easily bounded from above by
		\begin{align*}
			E_2^1&= L^2\sum_{i=1}^{\ell}\frac{K_h(\textbf{r}_k(x)-\textbf{r}_k(X_i))\|\textbf{r}_k(x)-\textbf{r}_k(X_i)\|^2}{\sum_{j=1}^{\ell}K_h(\textbf{r}_k(x)-\textbf{r}_k(X_j))}\mathds{1}_{\{D_h^{\ell}(x)>0\}}\mathds{1}_{\{\|\textbf{r}_k(x)-\textbf{r}_k(X_i)\|<R_Kh^{\beta}\}}\\
			&\leq L^2h^{2\beta}R_K^2\sum_{i=1}^{\ell}\frac{K_h(\textbf{r}_k(x)-\textbf{r}_k(X_i))}{\sum_{j=1}^{\ell}K_h(\textbf{r}_k(x)-\textbf{r}_k(X_j))}\mathds{1}_{\{D_h^{\ell}(x)>0\}}\\
			&= L^2h^{2\beta}R_K^2.
		\end{align*}
		Therefore, its expectation is simply bounded by the same upper bound i.e.,
		\begin{align}
			\label{eq:boundE21}
			\mathbb{E}(E_2^1)\leq L^2h^{2\beta}R_K^2
		\end{align}
		\item Secondly, we bound the second term $E_2^2$ using the tail assumption of the kernel $K$ given equation~\eqref{eq:assumption}, thus for any $h>0$:
		
		\begin{align*}
			E_2^2&= L^2\sum_{i=1}^{\ell}\frac{K_h(\textbf{r}_k(x)-\textbf{r}_k(X_i))\|\textbf{r}_k(x)-\textbf{r}_k(X_i)\|^2\mathds{1}_{\{D_h^{\ell}(x)>0\}}}{\sum_{j=1}^{\ell}K_h(\textbf{r}_k(x)-\textbf{r}_k(X_j))}\mathds{1}_{\{\|\textbf{r}_k(x)-\textbf{r}_k(X_i)\|\geq h^{\beta}R_K\}}\\
			&\leq L^2h^2\sum_{i=1}^{\ell}\frac{K_h(\textbf{r}_k(x)-\textbf{r}_k(X_i))\|(\textbf{r}_k(x)-\textbf{r}_k(X_i))/h\|^2\mathds{1}_{\{D_h^{\ell}(x)>0\}}}{\sum_{j=1}^{\ell}K_h(\textbf{r}_k(x)-\textbf{r}_k(X_j))}\times\\
			&\quad\ \mathds{1}_{\{(\|\textbf{r}_k(x)-\textbf{r}_k(X_i))/h\|\geq R_K/h^{1-\beta}\}}
		\end{align*}
		\begin{align*}
			&\leq \frac{h^2L^2}{b}\sum_{i=1}^{\ell}\frac{C_Ke^{-\|(\textbf{r}_k(x)-\textbf{r}_k(X_i))/h\|^{\alpha}}\|(\textbf{r}_k(x)-\textbf{r}_k(X_i))/h\|^2}{\sum_{j=1}^{\ell} \mathds{1}_{\{\|\textbf{r}_k(x)-\textbf{r}_k(X_j)\|<h\rho\}}}\times\\
			&\quad\ \mathds{1}_{\{\|(\textbf{r}_k(x)-\textbf{r}_k(X_i))/h\|\geq R_K/h^{1-\beta}\}}\mathds{1}_{\{C_h^{\ell}(x)>0\}}.
		\end{align*}
		As for any $\alpha>0$, $t\mapsto \lambda(t)=t^2e^{-t^{\alpha}}$ is strictly decreasing for all $t\geq (2/\alpha)^{1/\alpha}$. Thus, for $h>0$ small enough such that $R_K/h^{1-\beta}\geq(2/\alpha)^{1/\alpha}$, one has
		\begin{align*}
			E_2^2&\leq \frac{h^2L^2C_K}{b}\sum_{i=1}^{\ell}\frac{(R_K/h^{1-\beta})^2e^{-(R_K/h^{1-\beta})^{\alpha}}\mathds{1}_{\{\|(\textbf{r}_k(x)-\textbf{r}_k(X_i))/h\|\geq R_K/h^{1-\beta}\}}}{\sum_{j=1}^{\ell} \mathds{1}_{\{\|\textbf{r}_k(x)-\textbf{r}_k(X_j)\|<h\rho\}}}\mathds{1}_{\{C_h^{\ell}(x)>0\}}\\
			&\leq \frac{h^{2\beta}L^2C_KR_K^2e^{-R_K^{\alpha} h^{-\alpha(1-\beta)}}}{b}\sum_{i=1}^{\ell}\frac{\mathds{1}_{\{\sum_{j=1}^{\ell} \mathds{1}_{\{\|\textbf{r}_k(x)-\textbf{r}_k(X_j)\|<h\rho\}}>0\}}}{\sum_{j=1}^{\ell} \mathds{1}_{\{\|\textbf{r}_k(x)-\textbf{r}_k(X_j)\|<h\rho\}}}\\
			&\leq \frac{\ell h^{2\beta}L^2C_KR_K^2e^{-R_K^{\alpha} h^{-\alpha(1-\beta)}}}{b}\times\frac{\mathds{1}_{\{\sum_{j=1}^{\ell}\mathds{1}_{\{\|\textbf{r}_k(x)-\textbf{r}_k(X_j)\|<h\rho\}}>0\}}}{\sum_{j=1}^{\ell}\mathds{1}_{\{\|\textbf{r}_k(x)-\textbf{r}_k(X_j)\|<h\rho\}}}.
		\end{align*}
		Applying the result of inequality~(\ref{eq:bound1}), one has
		\begin{align}
			\mathbb{E}(E_2^2)&\leq\frac{\ell h^{2\beta}L^2C_KR_K^2e^{-R_K^{\alpha} h^{-\alpha(1-\beta)}}}{b}\times\frac{C_0}{h^M(\ell+1)}\nonumber\\
			&\leq C_1h^{2\beta-M}e^{-R_K^{\alpha}h^{-\alpha(1-\beta)}}\label{eq:boundE22}
		\end{align}
		for some $C_1>0$.
		\item Lastly with $A_h(x)$ defined in (\ref{eq:Axh}), we bound the expectation of $E_2^3$ by,
		\begin{align}
			\mathbb{E}(E_2^3)&\leq \mathbb{E}\Big[(g^*(\textbf{r}_k(x)))^2\mathds{1}_{\{C_h^{\ell}(x)=0\}}\Big]\nonumber\\
			&\leq \sup_{u\in\mathbb{R}^d}(g^*(\textbf{r}_k(u)))^2\mathbb{E}\Big[\mathds{1}_{\{C_h^{\ell}(x)=0\}}\Big]\nonumber\\
			&=\sup_{u\in\mathbb{R}^d}(g^*(\textbf{r}_k(u)))^2(1-\mu(A_h(x)))^{\ell}\nonumber\\
			&\leq \sup_{u\in\mathbb{R}^d}(g^*(\textbf{r}_k(u)))^2e^{-\ell\mu(A_h(x))}\nonumber\\
			&\leq \sup_{u\in\mathbb{R}^d}(g^*(\textbf{r}_k(u)))^2\frac{\ell\mu(A_h(x))e^{-\ell\mu(A_h(x))}}{\ell\mu(A_h(x))}\nonumber\\
			&\leq \sup_{u\in\mathbb{R}^d}(g^*(\textbf{r}_k(u)))^2\frac{\max_{u\in\mathbb{R}^d}ue^{-u}}{\ell\mu(A_h(x))}\nonumber
		\end{align}
		\begin{align}
			&\leq \sup_{u\in\mathbb{R}^d}(g^*(\textbf{r}_k(u)))^2\frac{e^{-1}}{\ell\mu(A_h(x))}\nonumber\\
			&\leq \frac{C_2}{\ell\mu(A_h(x)))}\label{eq:boundE23}
		\end{align}
		for some $C_2>0$.
	\end{itemize}
	From \eqref{eq:boundE1}, \eqref{eq:boundE21}, \eqref{eq:boundE22} and \eqref{eq:boundE23}, one has
	\begin{align*}
		\mathbb{E}[|g_n(\textbf{r}_k(X))-g^*(\textbf{r}_k(X))|^2]&\leq\int_{\mathbb{R}^d}\mathbb{E}[|g_n(\textbf{r}_k(x))-g^*(\textbf{r}_k(x))|^2]\mu(dx) \\
		&\leq \int_{\mathbb{R}^d}\mathbb{E}(E_1+E_2^1+E_2^2+E_2^3)\mu(dx)\\
		&\leq \int_{\mathbb{R}^d}\Big[\frac{4R^2}{b}\Big(\delta+\frac{C_0}{h^M(\ell+1)}\Big) + L^2h^{2\beta}R_K^2 \\
		&\quad+ C_1h^{2\beta-M}e^{-R_K^{\alpha}h^{-\alpha(1-\beta)}}  + \frac{C_2}{\ell\mu(A_h(x)))}\Big]\mu(dx).
	\end{align*}
	Therefore, by following the same procedure of proving inequality~(\ref{eq:bound1}), one has
	\begin{align*}
		&\quad\mathbb{E}[|g_n(\textbf{r}_k(X))-g^*(\textbf{r}_k(X))|^2]\\
		&\leq \frac{4R^2}{b}\Big(\delta+\frac{C_0}{h^M(\ell+1)}\Big) + L^2h^{2\beta}R_K^2+C_1h^{2\beta-M}e^{-R_K^{\alpha}h^{-\alpha(1-\beta)}} + \int_{\mathbb{R}^d}\frac{C_2\mu(dx)}{\ell\mu(A_h(x)))}\\
		&\leq \frac{4R^2}{b}\Big(\delta+\frac{C_0}{h^M(\ell+1)}\Big)+ L^2h^{2\beta}R_K^2+C_1h^{2\beta-M}e^{-R_K^{\alpha}h^{-\alpha(1-\beta)}}\\
		&\quad+\sum_{j\in J_{h,M}}\int_{\|\textbf{r}_k(x)-x_j\|<h\rho}\frac{C_2\mu(dx)}{\ell\mu(\{v\in\mathbb{R}^d:\|\textbf{r}_k(v)-\textbf{r}_k(x)\|<h\rho\})}\\
		&\leq \frac{4R^2}{b}\Big(\delta+\frac{C_0}{h^M(\ell+1)}\Big)+ L^2h^{2\beta}R_K^2+C_1h^{2\beta-M}e^{-R_K^{\alpha}h^{-\alpha(1-\beta)}}\\
		&\quad+\sum_{j\in J_{h,M}}\int_{\|\textbf{r}_k(x)-x_j\|<h\rho}\frac{C_2\mu(dx)}{\ell\mu(\{v\in\mathbb{R}^d:\|\textbf{r}_k(v)-x_j\|<h\rho\})}\\
		&\leq \frac{4R^2}{b}\Big(\delta+\frac{C_0}{h^M(\ell+1)}\Big)+ L^2h^{2\beta}R_K^2+C_1h^{2\beta-M}e^{-R_K^{\alpha}h^{-\alpha(1-\beta)}}\\
		&\quad+\frac{C_2}{\ell}\sum_{j\in J_{h,M}}\frac{\mu(\{v\in\mathbb{R}^d:\|\textbf{r}_k(v)-x_j\|<h\rho\})}{\mu(\{v\in\mathbb{R}^d:\|\textbf{r}_k(v)-x_j\|<h\rho\})}\\
		&\leq\frac{4R^2}{b}\Big(\delta+\frac{C_0}{h^M(\ell+1)}\Big)+ L^2h^{2\beta}R_K^2+C_1h^{2\beta-M}e^{-R_K^{\alpha}h^{-\alpha(1-\beta)}}+\frac{C_2|J_{h,M}|}{\ell}\\
		&\leq\frac{4R^2}{b}\Big(\delta+\frac{C_0}{h^M(\ell+1)}\Big)+ L^2R_K^2h^{2\beta}+C_1h^{2\beta-M}e^{-R_K^{\alpha}h^{-\alpha(1-\beta)}}+\frac{C_2'}{h^M\ell}
	\end{align*}
	where $|J_{h,M}|$ denotes the number of balls covering the ball $B$ (introduced in the proof of $A.2$) by the cover $\{B_M(x_j,h\rho):j=1,2,...\}$. Similarly, one has $|J_{h,M}|\leq \frac{C_0}{h^M}$ for some constant $C_0>0$ proportional to the volume of $B$. Since $\delta>0$ is chosen arbitrarily and the third term of the last inequality decreases exponentially fast when $h\to0$ for any $\beta\in(0,1)$, hence, it is negligible comparing to other terms. Finally, with the choice of $h\propto \ell^{-1/(M+2\beta)}$, one has
	\begin{equation*}
		\label{eq:final2}
		\mathbb{E}[|g_n(\textbf{r}_k(X))-g^*(\textbf{r}_k(X))|^2]\leq \frac{\tilde{C}_1}{h^M\ell}+\tilde{C}_2h^{2\beta}\leq C\ell^{-2\beta/(M+2\beta)}.
	\end{equation*}
	for some $C>0$ independent of $\ell$ and for any positive $\beta<1$ chosen arbitrarily. Thus, by letting $\beta\to1$, we obtain the desire result:
	\begin{equation*}
		\mathbb{E}[|g_n(\textbf{r}_k(X))-g^*(\textbf{r}_k(X))|^2]\leq C\ell^{-2/(M+2)}.
	\end{equation*}
\end{proofthm}
\QED

\section*{Acknowledgments}

The author gratefully acknowledges the support of Prof. Aurélie Fischer and Prof. Mathilde Mougeot for valuable feedback and suggestions during the process of writing this article.



\end{document}